\documentclass[reqno,11pt,a4paper]{amsart}

\usepackage{amssymb}

\numberwithin{equation}{section}
\newtheoremstyle{ttheorem}%
       {1.8ex\@plus1ex}                
       {2.1ex\@plus1ex\@minus.5ex}      
       {\itshape}           
       {0pt}                   
       {\bfseries}          
       {.}                  
       {.5em}               
       {}                

\newtheoremstyle{ddefinition}%
       {1.8ex\@plus1ex}                
       {2.1ex\@plus1ex\@minus.5ex}      
       {}           
       {0pt}                   
       {\bfseries}           
       {.}                  
       {.5em}               
       {}                

\newtheoremstyle{rremark}%
       {1.8ex\@plus1ex}                
       {2.1ex\@plus1ex\@minus.5ex}      
       {\normalfont}        
       {0pt}                   
       {\bfseries}           
       {.}                  
       {.5em}               
       {}                   

\theoremstyle{ttheorem}
\newtheorem{theorem}{Theorem}[section]
\newtheorem{lemma}[theorem]{Lemma}

\newtheorem{cor}[theorem]{Corollary}

\theoremstyle{ddefinition}
\newtheorem{definition}[theorem]{Definition}

\theoremstyle{rremark}
\newtheorem{remark}[theorem]{Remark}
\newtheorem{myremarks}[theorem]{Remarks}

\newenvironment{remarks}{\begin{myremarks}\begin{nummer}}%
    {\end{nummer}\end{myremarks}}

\newcounter{numcount}
\newcommand{\labelnummer}{\mbox{\normalfont (\roman{numcount})}}%

\makeatletter

\newenvironment{nummer}%
  {\let\curlabelspeicher\@currentlabel%
    \begin{list}{\labelnummer}%
      {\usecounter{numcount}\leftmargin0pt%
        \topsep0.5ex\partopsep2ex\parsep0pt\itemsep0ex\@plus1\p@%
        \labelwidth2.5em\itemindent3.5em\labelsep1em%
      }%
    \let\saveitem\item%
    \def\item{\saveitem%
      \def\@currentlabel{\curlabelspeicher$\,$\labelnummer}}%
    \let\savelabel\label%
    \def\label##1{\savelabel{##1}%
      \@bsphack%
        \ifmmode\else%
          \protected@write\@auxout{}%
          {\string\newlabel{##1item}{{\labelnummer}{\thepage}}}%
        \fi%
      \@esphack%
    }%
  }{\end{list}}%

\def\itemref#1{\expandafter\@setref\csname r@#1item\endcsname%
  \@firstoftwo{#1}}%

\def\section{\@startsection{section}{1}%
  \z@{1.3\linespacing\@plus\linespacing}{.5\linespacing}%
  {\normalfont\scshape\centering}}

\makeatletter

\def\one{\@ifundefined{comp}{\kern.5pt\leavevmode\hbox{\upshape{\small1\kern-3.35pt\normalsize1}}}{\kern.5pt\mathbb{1}}}%

\renewcommand{\le}{\leqslant}
\renewcommand{\ge}{\geqslant}
\renewcommand{\d}{\mathrm{d}}

\newcommand{\e}{\mathrm{e}}
\newcommand{\Chi}{\raisebox{.4ex}{$\chi$}}

\DeclareMathOperator{\supp}{\mathrm{supp}}
\DeclareMathOperator{\spec}{\mathrm{spec}}
\DeclareMathOperator{\specess}{\mathrm{spec}_{\mathrm{ess}}}
\DeclareMathOperator{\tr}{\mathrm{tr}\kern1pt}
\DeclareMathOperator{\dist}{\mathrm{dist}}
\DeclareMathOperator{\vol}{\mathrm{vol}}
\DeclareMathOperator{\diam}{\mathrm{diam}}

\newcommand{\cE}{\mathcal{E}}
\newcommand{\cG}{\mathcal{G}}
\newcommand{\hcG}{\widehat{\cG}}
\newcommand{\hXcG}{\widehat{X}_{\cG}}
\newcommand{\cV}{\mathcal{V}}

\renewcommand{\AA}{\mathbb{A}}
\newcommand{\PP}{\mathbb{P}}
\newcommand{\RR}{\mathbb{R}}
\newcommand{\QQ}{\mathbb{Q}}
\newcommand{\NN}{\mathbb{N}}
\newcommand{\ZZ}{\mathbb{Z}}

\def\per{.}
\def\HarvardComma{}

\newcounter{aucount}
\newif\ifedplural
\newif\ifper\pertrue

\def\au#1#2{{#1 #2}}
\def\lau#1#2{{#1 #2},}
\def\ed#1#2{\ifnum\theaucount=0\relax\fi{#1 #2}\addtocounter{aucount}{1}}
\def\led#1#2{\ifnum\theaucount=0\relax\edpluralfalse\else\edpluraltrue\fi{#1
    #2} (\editorname.),\setcounter{aucount}{0}}
\def\editorname{\ifedplural Eds\else Ed\fi}

\def\et{\ifnum\theaucount=1\else\HarvardComma\fi{} and\ }
\def\ti#1{#1,\ifper\fi\pertrue}

\def\bti{\@ifnextchar[\bbti\bbbti}
\def\bbti[#1]#2{\emph{#2}, #1,}
\def\bbbti#1{\emph{#1},}

\def\z{\@ifnextchar[\zz\zzz}
\def\zz[#1]#2#3#4#5{\perfalse\emph{#2} \textbf{#3}, #4 \ifx
  @#5@\relax\else (#5)\fi [#1]\ifper\per\fi\pertrue} 
\def\zzz#1#2#3#4{\emph{#1} \textbf{#2}, #3 \ifx @#4@\relax\else
  (#4)\fi\ifper\per\fi\pertrue}

\def\pub{\@ifstar\pubstar\pubnostar}
\def\pubnostar{\@ifnextchar[\@@pubnostar\@pubnostar}
\def\@@pubnostar[#1]#2#3#4{#2, #3, #4, #1\ifper\per\fi\pertrue}
\def\@pubnostar#1#2#3{#1, #2, #3\ifper\per\fi\pertrue}
\def\pubstar[#1]#2#3#4{\perfalse #2, #3, #4 [#1]\per\pertrue}

\makeatother

\begin{document}


\title[Random colourings of aperiodic graphs]{Random colourings of
  aperiodic graphs:\\ Ergodic and spectral
  properties}

\author[P.\ M\"uller]{Peter M\"uller}
\address{Institut f\"ur Theoretische Physik,
  Georg-August-Universit\"at G\"ottingen,
  Friedrich-Hund-Platz~1,
  37077 G\"ottingen, Germany}
\email{peter.mueller@physik.uni-goe.de}

\author[C.\ Richard]{Christoph Richard}        
\address{Fakult\"at f\"ur Mathematik,
  Universit\"at Bielefeld,
  Universit\"atsstr.\ 25,
  Postfach 10 01 31,
  33501 Bielefeld, Germany}
\email{richard@math.uni-bielefeld.de}


\begin{abstract}
  We study randomly coloured graphs embedded into Euclidean space,
  whose vertex sets are infinite, uniformly discrete subsets of finite
  local complexity. We construct the appropriate ergodic dynamical
  systems, explicitly characterise ergodic measures, and prove an
  ergodic theorem. For covariant operators of finite range defined on
  those graphs, we show the existence and self-averaging of the
  integrated density of states, as well as the non-randomness of the
  spectrum. Our main result establishes Lifshits tails at the lower
  spectral edge of the graph Laplacian on bond percolation subgraphs,
  for sufficiently small probabilities. Among other assumptions, its
  proof requires exponential decay of the cluster-size distribution
  for percolation on rather general graphs.
\end{abstract}

\maketitle


\section{Introduction}

Studying ensembles of random graphs is a broad subject with many
different facets. One of them, spectral properties of random graphs,
has found increasing interest in recent years. Its goal is to
determine spectral properties of the graph Laplacian, or of similar
operators associated with the graph, and to investigate their relation
to the graph structure. Erd\H{o}s-R\'enyi random graphs constitute one
class of examples, for which such types of results are known by now
\cite{KhSh04,KhKi06,BrDe06}. 

Another class of random graphs consists of those generated by a
percolation process from an underlying graph (``base graph''), which
is embedded into $d$-dimensional Euclidean space $\RR^{d}$. Standard
Bernoulli (bond- or site-) percolation subgraphs of the
$d$-dimensional hypercubic lattice are the prime example in this
category \cite{Gri99}. Here, ergodicity with respect to translations
has fundamental consequences, such as non-randomness of the spectrum,
as well as existence and self-averaging of the integrated density of
states \cite{Ves05, KiMu06}. The behaviour of the integrated density
of states near the edges of the spectrum requires a more detailed
understanding. Lifshits-tail behaviour was found in the
non-percolating phase \cite{KiMu06}, while the percolating cluster may
give rise to a van Hove asymptotics \cite{MuSt07}. In this context,
techniques from the theory of random Schr\"odinger operators have
turned out to be very efficient. Furthermore, the connection to the
theory of random walks in random environments \cite{Hug96,Bar04} was
exploited. Very recently, results of \cite{KiMu06,MuSt07} have been
extended to amenable Cayley graphs \cite{AnVe07a, AnVe07c}. There, it
is invariance under the appropriate group action, which replaces
translational invariance.

But how important is the automorphism group of the base graph for the
spectral asymptotics of its percolation subgraphs? To pursue this
question, we consider base graphs whose vertex sets are given by
infinite, uniformly discrete subsets of $\RR^{d}$, with the property
of finite local complexity (see Definition~\ref{FLC} below). Examples
include quasiperiodic tilings such as a Penrose tiling (see
\cite{BaMo00} for a recent monograph on quasiperiodic point sets),
more generally, tilings with a finite set of prototiles \cite{GrSh89},
but also random tiling ensembles \cite{RHHB98}.  Typically, none of
these enjoys invariance under an appropriate group action. Ergodic and
spectral properties of the base graphs were first derived by
\cite{Hof93, Hof95}, and significantly extended by \cite{LeSt03,
  KlLe03, LeSt05}, using methods from dynamical systems. In this
paper, we supply these base graphs with a random colouring and study
their spectral properties.  The main result of this paper,
Theorem~\ref{main}, goes beyond basic ergodic spectral properties and
establishes Lifshits-tail behaviour at the lower spectral edge for the
graph Laplacian on percolation subgraphs.

Our proof of this result involves three preparatory steps, each of
which is interesting in its own. The first step belongs to the realm
of dynamical systems theory, the second to spectral theory,
and the third to percolation theory.

\emph{(i) \; Construct the appropriate ergodic dynamical systems,
  explicitly characterise ergodic measures and prove an ergodic
  theorem.}  Given an ergodic measure on the dynamical system of the
base graphs, we will explicitly construct an ergodic measure for
corresponding randomly coloured graphs, following ideas of
\cite{Hof98}.  The main result of this step is an ergodic theorem
(Theorem~\ref{theo:perg}) for dynamical systems associated with randomly
coloured graphs. It extends \cite{Hof98}, where colourings of
aperiodic Delone graphs with strictly ergodic dynamical system have
been studied. Our setting covers the full range from periodic
structures to random tilings. Moreover, we do not require relative
denseness of the vertex sets, thereby including examples such
as the visible lattice points \cite{BaMoPl00} in our setup. 
Apparently, some of the technical problems we had to overcome are 
closely related to ones in \cite{BaaZin07}, where diffraction properties 
of certain random point sets, including percolation subsets, have been 
investigated very recently.

\emph{(ii) \; Derive ergodic spectral properties of covariant,
  finite-range operators on randomly coloured aperiodic graphs.}
Theorem~\ref{maclimit} characterises the integrated density of states
of such an operator by a macroscopic limit.  Theorem~\ref{spec-ids}
states the non-randomness of the spectrum of the operator and relates
it to the set of growth points of the integrated density of states.
In particular, the theorems guarantee that there are no exceptional
instances to their statements for uniquely ergodic systems.  We
provide elementary proofs of Theorems~\ref{maclimit} and
\ref{spec-ids}. In the absence of a colouring, corresponding results
have been derived in \cite{Hof93,Hof95, LeSt03,LeSt05}, mainly in the
strictly ergodic or in the uniquely ergodic case. 

\emph{(iii) \; Establish exponential decay of the cluster-size
  distribution in the non-percolating phase for general graphs.}  We
derive an elementary exponential-decay estimate for the probability to
find an open path from the centre to the complement of a large ball.
Unfortunately, this estimate holds only for sufficiently small bond 
probabilities. For these probabilities, the decay of the cluster-size 
distribution then follows from that estimate, by verifying that the 
corresponding arguments in \cite{Gri99} apply also in our general 
setting. Exponential decay throughout the non-percolating phase for
quasi-transitive graphs has been proved recently \cite{AnVe07b}.
Within our more general setup, an extension to higher bond
probabilities up to criticality remains a challenging open question,
see also the discussion in \cite{Hof98}.

The manuscript \cite{LeVe07}, which was finalised at the same time as
ours, establishes uniform convergence in the energy of the
finite-volume approximants to the integrated density of states under
rather general conditions. In particular, it applies to percolation on
Delone dynamical systems and thus improves on Theorem~\ref{maclimit}
under slightly different conditions. However, the validity of our
general Ergodic Theorem~\ref{theo:perg} is an open question in the
approach of \cite{LeVe07}. Using uniform convergence would not allow
to strengthen our main result on Lifshits tails in Theorem~\ref{main}.

Our paper is organised as follows. Section~\ref{secdyn} sets the
notation and introduces dynamical systems associated with uncoloured
graphs. This is a slight extension of the setup for Delone dynamical
systems, such as in \cite{LeMo02}.  In Section~\ref{sec:ran-col}, we
construct the dynamical systems for the corresponding randomly
coloured graphs and deal with step~(i).  Section~\ref{sec:operators}
introduces covariant operators of finite range on randomly coloured
graphs and treats step~(ii). Section~\ref{secLif} is devoted to our
main result on Lifshits tails together with its proof.
Section~\ref{sec:spectrum} contains the proof of
Theorem~\ref{spec-ids}, and Section~\ref{sec:percest} deals with step
(iii).

%
\section{Dynamical systems for graphs}
\label{secdyn}
%

For the basic notions involving graphs, we refer, for example, to the
textbook \cite{Die05}. We consider (simple) graphs $G= (\cV,\cE)$,
whose vertex sets $\cV\equiv \cV_{G}$ are countable subsets of
$\mathbb R^d$.  We say that $\cV$ is \emph{uniformly discrete} of
radius $r \in ]0,\infty[$, if any open ball of radius $r$ in $\RR^{d}$
contains at most one element of $\cV$. The vertex set is called
\emph{relatively dense} if there exists $R \in[0,\infty[$ such that
every closed ball of radius $R$ contains at least one vertex. The vertex set
is called a \emph{Delone set}, if it is both uniformly discrete and
relatively dense. The edge set $\cE\equiv\cE_{G}$ of $G$ is a subset
of the set of all unordered pairs of vertices. We denote an edge by
$e\equiv\{v,w\}$, where $v,w \in\cV$ with $v\neq w$. In other words,
we do not allow self-loops, nor multiple edges between the same pair
of vertices.

Recall that a graph $G'=(\mathcal{V}',\mathcal{E}')$ is called a
subgraph of $G$, in symbols $G'\subseteq G$, if $\cV'\subseteq \cV$
and $\cE'\subseteq \cE$. For $x\in\mathbb R^d$, the \emph{translated
  graph} $x+G$ has vertex set $x+\cV :=\{x+v:v\in\cV\}$ and edge set
$x+\cE := \bigl\{ \{x+v,x+w\}: \{v,w\} \in \cE \bigr\}$. Given any
Borel set $B\subseteq\RR^{d}$, the \emph{restriction} $G\wedge B$ of
$G$ to $B$ is the induced subgraph of $G$ with vertex set $\cV \cap
B$, that is, $\{u,v\}$ belongs to the edge set of $G\wedge B$, if and
only if $\{u,v\} \in \cE$ and $u,v\in \cV \cap B$.  If $B$ is bounded,
then $G\wedge B$ is called a \emph{$B$-pattern} (or simply a pattern)
of $G$. Two patterns $P,Q$ are called \emph{equivalent}, if $x+P=Q$
for some $x\in\mathbb R^d$. An \emph{$r$-pattern} is a pattern
$G\wedge B_r(v)$ for some $v\in \cV_G$. Here we have used the notation
$B_r(x)$ for the open ball of radius $r>0$ around $x\in
\mathbb{R}^{d}$. In particular, we set $B_r:=B_r(0)$. We write $|M|$
for the cardinality of a set $M$.

\begin{definition}
  \begin{nummer}
  \item
    \label{FLC}    
    A set $\cG$ of graphs is said to have \emph{finite local
      complexity}, if for every $r>0$
    \begin{equation}
      \big| \{(-v+G)\wedge B_r: v\in\cV_{G}, G\in\cG\} \big| <\infty.
    \end{equation}
    In particular, a single graph $G$ has finite local complexity if,
    for any given $r>0$, the number of its non-equivalent $r$-patterns
    is finite.
  \item 
    Let $G$ be a finite graph and $P \subseteq G$ a pattern of $G$. The
    \emph{number of occurrences}
    \begin{equation}
      \nu(P|G) := |\{x\in\RR^d : x+P \subseteq G\}|
    \end{equation}
    of $P$ in $G$ is the (finite) number of translates of $P$ in $G$.    
  \end{nummer}
\end{definition}

Geometric properties of some set of graphs $\cG$ are reflected by
properties of an associated dynamical system. This we introduce along
the lines of \cite{LeMo02}, where the case of Delone multi-sets was
considered. The statements of this section are proved by slight
adaptations of the arguments laid down in \cite{LeMo02,RaWo92,Sch00}. In
fact, examples of our setup include the Delone multi-sets of
\cite{LeMo02}, in which case $\cG$ is finite, vertex sets are Delone
sets and edge sets are empty. 

For simplicity, let us assume now that the vertex set of each
$G\in\cG$ is uniformly discrete. Following \cite{LeMo02,Hof95}, we
define a \emph{metric} on $\cG$ by setting
\begin{align}
  \label{metric}
  \mathrm{dist}(G,G') :=\min\Big\{ 2^{-1/2}, \inf\big\{ & \varepsilon>0:
  \text{~there exists~} x,y\in 
  B_\varepsilon: \nonumber\\ 
  & (x+G) \wedge B_{1/\varepsilon}=(y+G')\wedge
  B_{1/\varepsilon}\big\} \Big\}
\end{align}
for all $G,G' \in\cG$. In essence, two graphs are close, if they agree,
up to a small translation, on a large ball around the origin. 
Symmetry and the triangle inequality of the metric are seen to hold
for any set of graphs $\cG$. The  uniform discreteness assumption
ensures positive definiteness, but it is much stronger than what is
required. In fact, it would have been sufficient to
assume merely closedness of the vertex sets and of the edge sets
(with respect to a suitable metric on $\cE$). Now we define the
complete metric space
\begin{equation}
  X_{\cG} := \overline{ \{x +G: x \in\RR^{d}, G \in\cG \} },
\end{equation}
of all translates of graphs in $\cG$, where the metric \eqref{metric}
is used for completion. Later we will need to know that certain
properties of graphs in $\cG$ do not get lost in the closure.

\begin{lemma}
  Let $\cG$ be a set of graphs with uniformly discrete vertex sets.
  \begin{nummer}
    \item
      \label{dmax}
      If for some $d_{\mathrm{max}} \in \mathbb{N}$ the estimate
      \begin{equation}
        \sup_{v \in\cV_{G}} d_{G}(v) \le d_{\mathrm{max}}
      \end{equation}
      holds for all $G\in \cG$, then it holds also for all
      $G \in X_{\cG}$.
    \item 
      \label{notfinite}
      In addition to uniform discreteness, assume there is some finite
      $R >0$, such that the vertex sets of all $G\in\cG$ 
      are relatively dense with radius $R$, and that
      \begin{equation}
        \ell_{\mathrm{max}} := \sup \bigl\{|u-v|: G\in\cG, \{u,v\} \in
        \cE_{G}\bigr\} < \infty,  
      \end{equation}
      i.e., there exists a finite maximum bond length. Then, every $G\in
      X_{\cG}$ is infinite, and if no $G\in\cG$
      possesses a finite cluster, then no $G\in X_{\cG}$ possesses a
      finite cluster.
  \end{nummer}  
\end{lemma}

\begin{proof}
  By contradiction. 
  \begin{nummer}
  \item Assume there exists $G \in X_{\cG}$ and $v \in \cV_{G}$ such
    that $d_{G}(v) > d_{\mathrm{max}}$. Choose $0< \varepsilon <1/3$
    small enough, such that all neighbours of $v$ lie in the ball
    $B_{\varepsilon^{-1}-1}(v)$. Hence we have $d_{G}(v) = d_{G \wedge
      B_{1/\varepsilon}(v)}(v)$. Since $G$ is in the closure of $\cG$, there
    exists $G' \in\cG$ with $\dist(G, G') < \varepsilon$. Altogether, this
    implies $d_{G'\wedge B_{1/\varepsilon}(v)}(v)  >
    d_{\mathrm{max}}$, a contradiction.
  \item 
    The assumption of relative denseness implies that $X_{\cG}$
    contains only infinite graphs. So let us assume there exists $G
    \in X_{\cG}$ with a finite cluster $C$. Let $r_{C} \in ]0,\infty[$
    big enough such that $\cV_{C} \subset B_{r_{C}}$, and choose
    $\varepsilon> 0$ so small that $1/\varepsilon > r_{C} +
    \ell_{\mathrm{max}} +3$.  Since $G$ is in the closure of $\cG$, there
    exists $G' \in\cG$ with $\dist(G, G') < \varepsilon$, and $G'\wedge
    B_{1/\varepsilon}$ has a finite cluster that cannot merge with
    other clusters when removing the restriction to the ball
    $B_{1/\varepsilon}$.  
    \qed
  \end{nummer}  
  \def\qed{}
\end{proof}

Standard arguments show that the translation
group $\RR^{d}$ acts continuously on $X_{\cG}$, that is, the map $G
\mapsto x+ G$ is continuous for every $x\in\RR^{d}$. Thus, the triple
$(X_{\cG}, \mathbb R^d,+)$ constitutes a topological dynamical system.
The following result can be proved along the lines of
\cite{RaWo92,Sch00}.

\begin{lemma}
  \label{lem:flc}
  Let $\cG$ be set of graphs with uniformly discrete vertex
  sets. Then, $X_{\cG}$ is compact if and only if $\cG$ has finite
  local complexity. \qed
\end{lemma}

As above, uniform discreteness can be replaced by closedness of all
vertex and edge sets, without jeopardising the validity of
Lemma~\ref{lem:flc}.  But from now on, we assume that uniform
discreteness holds even uniformly in $\cG$, that is, there exists
$r>0$ such that $\cV_{G}$ is uniformly discrete of radius $r$, for all
$G\in\cG$.

Compactness of $X_{\cG}$ implies the existence of ergodic probability
measures on the Borel-sigma algebra of $X_{\cG}$, in other words $X_{\cG}$
is \emph{ergodic} (w.r.t.\ translations). Recall that a topological
dynamical system is called \emph{uniquely ergodic}, if it carries exactly
one ergodic measure. Ergodic theorems for a compact
dynamical system with $\mathbb R^d$-action are given in
\cite[Thms.~4.2 and~2.6]{LeMo02}, see also \cite[Thm.~1]{LeSt05} for a
stronger statement in the case of minimal ergodic systems.
We quote a version patterned after \cite{LeMo02} and introduce 
\emph{cylinder sets}
\begin{equation}
  \label{cyl-set}
  \Xi_{P, U} := \{G \in X_{\cG}: x+P \subseteq G \text{~for some~}
  x\in U\} \subset X_{\cG}.
\end{equation}
Here, $P$ is a pattern of some graph $G\in\cG$,
and $U\subseteq\RR^{d}$ is a Borel set. We write $\vol(U)$ for its
Lebesgue measure. 

\begin{theorem}
  \label{basic-erg}
  Let $\cG$ be a set of graphs of finite local complexity and with
  uniformly discrete vertex sets of radius $r>0$. Fix an ergodic
  probability measure $\mu$ on $X_{\cG}$. Then, given any function
  $\phi\in \mathrm{L}^1(X_{\cG},\mu)$, the limit
  \begin{equation}
    \label{basic-erg-limit}
    \lim_{n\to\infty} \frac{1}{\vol(B_{n})}  \int_{B_{n}}\!\d x\;
    \phi(x+G)  
    =   \int_{X_{\cG}}
    \!\d\mu(F) \; \phi (F)
  \end{equation}
  exist for $\mu$-a.a.\ $G \in X_{\cG}$.  If $X_{\cG}$ is even
  uniquely ergodic and, in addition, if $\phi$ is either continuous or
  a linear combination of indicator functions of cylinder sets, then
  the limit \eqref{basic-erg-limit} exists for \emph{all} $G\in
  X_{\cG}$. \qed
\end{theorem}

The ergodic theorem can be used to analyse the vertex density of 
graphs and the asymptotic number of
occurrences of patterns in graphs.

\begin{cor}
  \label{sammelcor}
  Let $\cG$ be a set of graphs of finite local complexity and with
  uniformly discrete vertex sets of radius $r>0$.  Fix an ergodic
  probability measure $\mu$ on $X_{\cG}$. 
  \begin{nummer}
  \item 
    \label{density}
    Then there is a Borel set $\overline{X} \subseteq X_{\cG}$ of full
    $\mu$-measure, $\mu(\overline{X})=1$, such that the \emph{vertex
      density}
    \begin{equation}
      \label{rho-def}
      \varrho := \lim_{n\to\infty} \frac{|\mathcal{V}_{G} \cap
        B_{n}|}{\vol(B_{n})}  = \int_{X_{\cG}}\!\d\mu(F) \,
      \sum_{v\in\mathcal{V}_{F}}  \psi(v) 
    \end{equation}
    exists for all $G \in \overline{X}$. Here, $\psi: \RR^d\to\RR_{\ge 0}$
    is any continuous function with support in $B_{r}$ and
    $\int_{\RR^{d}}\d x\, \psi(x) =1$ (``mollifier''). In particular,
    $\varrho\in [0,1/\vol(B_{r})]$ is independent of $G\in \overline{X}$
    and of the choice of the mollifier. If $X_{\cG}$ is
    even uniquely ergodic, then the above statements hold with
    $\overline{X}=X_{\cG}$.
  \item 
    \label{rho-infty}
    The statements of Part~\itemref{density} apply also to the
    \emph{density of vertices belonging to infinite components}
    \begin{equation}
      \label{rho-infty-def}
      \varrho_{\infty} := \lim_{n\to\infty} \frac{|\mathcal{V}_{G,\infty} \cap
        B_{n}|}{\vol(B_{n})}  = \int_{X_{\cG}}\!\d\mu(F) \,
      \sum_{v\in\mathcal{V}_{F,\infty}}  \psi(v) .
    \end{equation}
    Here, $\cV_{G,\infty} := \{v\in\cV_{G} : |C_{v}| = \infty\}$, with
    $C_{v}$ denoting the cluster of $G$ which $v\in\cV_{G}$ belongs to.
  \item 
    \label{lem:freq}
    Let $P$ be a pattern of some graph $G\in\cG$.  Then there is a
    Borel set $\overline{X}\subseteq X_{\cG}$ of full $\mu$-measure, such that
    the \emph{pattern frequency}
    \begin{equation}
      \nu(P) := \lim_{n\to\infty} \frac{\nu(P|G\wedge
        B_{n})}{\vol(B_{n})}   
    \end{equation} 
    exists for all $G \in \overline{X}$ and is independent of $G$. 
    The dynamical system $X_{\cG}$ is uniquely ergodic, if and only if,
    given any pattern $P$ of a graph in $\cG$, the limit  
    \begin{equation}
      \label{upf-eq}
      \nu(P) := \lim_{n\to\infty}\frac{\nu(P|G\wedge B_{n}(a))}
      {\vol(B_{n})}
    \end{equation} 
    exists uniformly in $G \in X_{\cG}$ and in $a\in\RR^d$, and is
    independent of $G$ and $a$.
  \item 
    \label{lem:cylinder}
    Let $P$ be a pattern of some graph $G\in\cG$, and let
    $U\subset\RR^d$ be a Borel set with diameter
    $\mathrm{diam}(U)<r$. Then, the probability of the associated
    cylinder set \eqref{cyl-set} is given by
    \begin{equation}
      \mu(\Xi_{P, U})= \vol(U) \, \nu(P).
    \end{equation}
  \end{nummer}
\end{cor}

\begin{remark}
  \begin{nummer}
  \item 
    \label{dens-well-def}
    The sum over $v$ in the $\mu$-integral in \eqref{rho-def} contains
    at most one term, because the mollifier $\psi$ is supported in the
    ball $B_{r}$, where $r$ is the radius of uniform discreteness of
    the graphs.
  \item 
    The criterion \eqref{upf-eq} for unique ergodicity in
    Lemma~\ref{lem:freq} is often referred to as \emph{uniform pattern 
      frequencies}.
  \item 
    \label{rhopos}
    For later reference we give two different conditions that imply
    $\varrho_{\infty} >0$: \quad(1)\quad The situation described in
    Lemma~\ref{notfinite}, assuming that no $G\in\cG$ possesses a
    finite cluster. \quad(2)\quad The dynamical system $X_{\cG}$ is uniquely
    ergodic and there exists $G \in X_{\cG}$ such that
    $|\mathcal{V}_{G,\infty} \cap B_{n}| / \vol(B_{n})$ has a strictly
    positive limit as $n\to\infty$.
  \end{nummer}  
\end{remark}

\begin{proof}[Proof of Corollary~\ref{sammelcor}.]
  Part~\itemref{density} of the corollary follows from an application
  of Theorem~\ref{basic-erg} to the continuous function $\phi(G) :=
  \sum_{v\in\cV_{G}} \psi(v)$ and the relation
  \begin{equation}
    |\cV_{G} \cap B_{n}| = \int_{B_{n}}\!\d x\; \phi(x+G) +
    \mathcal{O}(n^{d-1}). 
  \end{equation}
  The latter reveals the independence of the right-hand side of
  \eqref{rho-def} on the particular choice of the mollifier $\psi$. For
  Part~\itemref{rho-infty}, one has to replace $\cV_{G}$ by
  $\cV_{G,\infty}$ in the argument.

  Parts~\itemref{lem:freq} and~\itemref{lem:cylinder} follow from repeating
  the arguments in \cite[Lemma~4.3 and Thm.~2.7]{LeMo02}, where the
  case of Delone multi-sets was treated. 
\end{proof}

\begin{definition}
  \label{posfreq}
  Let $\cG$ be a set of graphs of finite local complexity and with
  uniformly discrete vertex sets of radius $r>0$. We say that the
  dynamical system $X_{\cG}$ satisfies the \emph{positive lower frequency
    condition}, if for every  $G\in X_{\cG}$ and every pattern
  $P\subset G$ 
  \begin{equation}
    \label{posfreq-eq}
    \liminf_{n\to\infty} \;\frac{\nu(P|G\wedge B_{n})}{\vol(B_{n})} > 0.
  \end{equation}
  Loosely spoken, any pattern $P$ that occurs once in $G$, does so
  sufficiently often.
\end{definition}

\begin{remarks}  \label{mini-posfreq}
\item 
  Minimality of $X_{\cG}$, which is equivalent to repetitivity for
  Delone systems of finite local complexity \cite[Thm.~3.2]{LaPl03}, implies
  that the positive lower frequency condition holds.
\item 
  If, in addition to the positive lower frequency condition, one assumes
  that $X_{\cG}$ is uniquely ergodic, then, in view of
  Corollary~\ref{lem:freq}, the $\liminf$ in \eqref{posfreq-eq} equals the
  pattern frequency $\nu{(P)}$. Moreover, in this case the system is
  minimal and thus strictly ergodic, compare \cite{LaPl03}.
\end{remarks}

%
\section{Ergodic properties of randomly coloured graphs}
\label{sec:ran-col}
%

In this section, we supply the graphs of the previous section with a
random colouring, and derive a corresponding extension of the Ergodic
Theorem~\ref{basic-erg}.

We fix a finite, nonempty set $\AA$, equipped with the discrete topology,
which we call the set of available \emph{colours}. For definiteness, we
consider only random edge colourings of graphs. But all results of this and
the next section remain valid in the case of a random colouring of vertices,
and of a random colouring of both edges and vertices. This is merely
a matter of notation.

For a given a graph $G$, we define the probability space
$\Omega_{G}:=\mbox{\Large$\times$}_{e\in\cE_{G}}\AA$, equipped with the
$|\cE_{G}|$-fold product sigma-algebra $\bigotimes_{e\in\cE_{G}}\, 2^{\AA}$ of
the power set of $\AA$ and the product probability measure $\PP_{G} :=
\bigotimes_{e \in \cE_{G}} \PP_0$. Here, $\PP_{0}$ is some fixed probability
measure on $\AA$.  In other words, colours are distributed identically and
independently to all edges, and the elementary event $\omega \equiv
(\omega_{e})_{e \in \cE_{G}}\in\Omega_{G}$ specifies a particular realisation
of colours assigned to the edges of $G$.

At first we are going to extend the framework of the previous section to
coloured graphs. Given a graph $G$ and $\omega\in\Omega_{G}$, the pair
$G^{(\omega)} \equiv (G,\omega)$ is called a \emph{coloured graph}. For any
Borel set $B\subseteq \RR^{d}$, we define the restriction $\omega\wedge B \in
\mbox{\Large$\times$}_{e\in\cE_{G\wedge B}}\AA$ as the image of $\omega$ under
the canonical projection from $\Omega_{G}$ to
$\mbox{\Large$\times$}_{e\in\cE_{G\wedge B}}\AA$. Likewise, given any
$x\in\RR^{d}$, the translated colour realisation $x+\omega \in \Omega_{x+G}$
is defined component-wise by $(x+\omega)_{x+e} := \omega_{e}$ for all
$e\in\cE_{G}$. We define the translation and restriction of a coloured
graph in the natural way
\begin{equation}
  x + G^{(\omega)} :=  (x+G)^{(x+\omega)}, \qquad\quad
  G^{(\omega)} \wedge B := (G\wedge B)^{(\omega\wedge B)},
\end{equation}
by shifting and truncating $\omega$ along with $G$. We write
$P^{(\eta)} \subseteq G^{(\omega)}$, if $P \subseteq G$ and $\eta_{e} =
\omega_{e}$ for all edges $e\in \cE_{P}$.
The notions of a pattern of a coloured graph and of the number of
occurrences of a finite coloured graph $P^{(\eta)}$ in $G^{(\omega)}$
translate accordingly from those in the previous section. 

For a given set of graphs $\cG$, we consider the induced set of
coloured graphs
\begin{equation}
 \hcG := \{G^{(\omega)}: \omega\in\Omega_{G}, G\in\cG\}. 
\end{equation}

\begin{remark}
  \label{transfer}
  Since $\AA$ provides only finitely many different colours,
  it follows that $\cG$ has finite local complexity, if and only if
  $\hcG$ has finite local complexity, that is if and only if 
  \begin{equation}
    |\{(-x+G^{(\omega)}) \wedge B_r:  x\in\cV_{G},
    G^{(\omega)}\in\hcG  \}|  <\infty
  \end{equation}
  for every $r>0$.   
\end{remark}

Replacing $G$ and $G'$ in the metric \eqref{metric} by
elements of $\hcG$, we obtain a metric on
$\hcG$. This metric is used in the completion of the metric space 
\begin{equation}
  \hXcG := \overline{ \{x +G^{(\omega)}: x \in\RR^{d},
    G^{(\omega)} \in\hcG \} \rule{0pt}{2.2ex}}.
\end{equation}
An alternative description of the space $\hXcG$ is provided by

\begin{lemma}
  Let $\cG$ be a set of graphs with uniformly discrete vertex sets.
  Then
  \begin{equation}
    \label{hat-equal}
    \hXcG   = \{ G^{(\omega)}: \omega\in\Omega_{G}, G \in X_{\cG}\}.
  \end{equation}
\end{lemma}

\begin{proof}
  To show the inclusion ``$\,\subseteq\,$'', it suffices to prove that the
  limit $\widehat{G}$ of an arbitrary convergent sequence from $\hcG$ is of
  the form $G^{(\omega)}$ for some $G\in X_{\cG}$ and some
  $\omega\in\Omega_{G}$. So assume that for every $\varepsilon >0$ there exist
  $x_{\varepsilon} \in\RR^{d}$, $y_{\varepsilon} \in B_{\varepsilon}$ and
  $G_{\varepsilon}^{(\omega_{\varepsilon})} \in \hcG$ such that $
  (x_{\varepsilon} +G_\varepsilon^{(\omega_{\varepsilon})}) \wedge
  B_{1/\varepsilon} = (y_{\varepsilon} + \widehat{G}) \wedge
  B_{1/\varepsilon}$.  Clearly, convergence of a sequence of coloured graphs
  implies convergence of the underlying uncoloured graphs, that is, there
  exists $G \in X_{\cG}$ such that $(x_{\varepsilon} +G_\varepsilon) \wedge
  B_{1/\varepsilon} = (y_{\varepsilon} + G) \wedge B_{1/\varepsilon}$ for all
  $\varepsilon >0$. We define $\omega\in\Omega_{G}$ as follows: given any
  $\varepsilon>0$, we set $\omega_{e} := \omega_{\varepsilon, y_{\varepsilon}
    -x_{\varepsilon} +e}$ for all $e \in\cE_{G\wedge
    B_{1/\varepsilon}(-y_{\varepsilon})}$. This choice is consistent in the
  sense that if $0<\varepsilon' < \varepsilon$, then $\omega_{\varepsilon,
    y_{\varepsilon} -x_{\varepsilon} +e} = \omega_{\varepsilon',
    y_{\varepsilon'} -x_{\varepsilon'} +e}$ for all $e\in \cE_{G\wedge
    B_{1/\varepsilon}(-y_{\varepsilon})}$. By choosing $\varepsilon$
  arbitrarily small, we obtain $\omega_{e}$ for all $e\in\cE_{G}$. It follows
  that $\widehat{G}=G^{(\omega)}$.

  To show the inclusion ``$\,\supseteq\,$'', consider $G^{(\omega)}$ for an
  arbitrary $G\in X_{\cG}$ and arbitrary $\omega \in\Omega_{G}$. Then, for
  every $\varepsilon >0$ there exist $x_{\varepsilon} \in\RR^{d}$,
  $y_{\varepsilon} \in B_{\varepsilon}$ and $G_{\varepsilon} \in \cG$ such
  that $ (x_{\varepsilon} +G_\varepsilon) \wedge B_{1/\varepsilon} =
  (y_{\varepsilon} + G) \wedge B_{1/\varepsilon}$.  Define
  $\omega_{\varepsilon,e}:=\omega_{x_\varepsilon -y_{\varepsilon} +e}$ for all
  $e \in\cE_{G_{\varepsilon}\wedge B_{1/\varepsilon}(-x_{\varepsilon})}$ and
  set $\omega_{\varepsilon, e}$ to an arbitrary value for the remaining edges.
  We then have $G_{\varepsilon}^{(\omega_{\varepsilon})} \in \hcG$ for all
  $\varepsilon >0$ and $x_{\varepsilon} +
  G_\varepsilon^{(\omega_\varepsilon)}\to G^{(\omega)}$ as
  $\varepsilon\downarrow 0$.
\end{proof}

There is an analogue to Lemma~\ref{lem:flc} in the previous
section. Recalling Remark~\ref{transfer}, it can be
formulated as

\begin{lemma}
  \label{lem:flc-hat}
  Let $\cG$ be a set of graphs with uniformly discrete vertex sets.
  Then, $\hXcG$ is compact if and only if $\cG$ has
  finite local complexity. \qed
\end{lemma}

The main goal of this section is to express an ergodic probability measure
$\widehat{\mu}$ on a compact space $\hXcG$ in terms of an ergodic probability
measure $\mu$ on $X_{\cG}$ and the probability measures $\mathbb P_G$. This
will be achieved in Theorem~\ref{theo:perg} at the end of this section.

As a preparation for Theorem~\ref{theo:perg}, we define cylinder sets of
$\hXcG$ in analogy those of $X_{\cG}$ in the previous section. Given a pattern
$P^{(\eta)}$ of some coloured graph in $\hcG$ and a Borel set
$U\subseteq\RR^{d}$, we set
\begin{equation}
  \widehat{\Xi}_{P^{(\eta)},U} := \{ G^{(\omega)} \in \hXcG :  x + P^{(\eta)}
  \subseteq G^{(\omega)} \text{~for some~} x\in U\}.
\end{equation}
The basic step in the construction of an ergodic measure on
$\hXcG$ is given by the following lemma, which extends
\cite[Lemma~3.1]{Hof98} to graphs which are not
necessarily aperiodic -- and also to more general measures $\PP_{G}$.
We employ the notation $\Chi_{S}$ for the indicator function
of some set $S$.

\begin{lemma}
  \label{lem:hoflemma}
  Let $\cG$ be a set of graphs of finite local complexity and with
  uniformly discrete vertex sets of radius $r>0$. Let $\mu$ be an
  ergodic probability measure on $X_{\cG}$. Then, given any pattern
  $P^{(\eta)}$ of some coloured graph in $\hcG$ and any Borel set
  $U\subseteq\RR^{d}$ with diameter $\diam(U) <r$, the limit
  \begin{equation}
    \label{erg-conv}
    \lim_{n\to\infty} \frac{1}{\vol(B_{n})} \int_{B_n}\!\d x\;
    \Chi_{\widehat{\Xi}_{P^{(\eta)},U}}(x+G^{(\omega)})
    =  \mu(\Xi_{P, U}) \, \PP_P(\eta)
  \end{equation}
  exists for $\mu$-a.a.\ $G \in X_{\cG}$ and for $\mathbb P_G$-a.a.\
  $\omega\in\Omega_G$. If $\hXcG$ is uniquely ergodic, then the limit
  \eqref{erg-conv} exists for \emph{all} $G\in X_{\cG}$ and for
  $\mathbb P_G$-a.a.\ $\omega\in\Omega_G$.
\end{lemma}

\begin{proof}
  It follows from the definition of cylinder sets that
  \begin{equation}
    \label{hof-start}
    \int_{B_n}\!\d x\;
    \Chi_{\widehat{\Xi}_{P^{(\eta)},U}}(x+G^{(\omega)})
    =  \vol(U) \; \nu(P^{(\eta)}|G^{(\omega)} \wedge B_{n})
    +\mathcal{O}(n^{d-1})
  \end{equation}
  asymptotically as $n\to\infty$, see also the proof of \cite[Thm.~2.7]{LeMo02}
  for the argument. We analyse the asymptotic behaviour of the expression
  \begin{equation}
    \label{count}
    \frac{\nu(P^{(\eta)}|G^{(\omega)} \wedge B_{n})}{\vol(B_{n})}\le
    \frac{\nu(P|G \wedge B_{n})}{\vol(B_{n})}.
  \end{equation}
  Ergodicity of $X_{\cG}$ implies
  \begin{equation}
    \label{freq-conv}
    \lim_{n\to\infty}  \frac{\nu(P|G \wedge B_{n})}{\vol(B_{n})} = \nu(P) =
  \frac{\mu(\Xi_{P,U})}{\vol(U)},
  \end{equation}
  for $\mu$-a.a.\ $G\in X_{\cG}$, resp.\ for all $G \in X_{\cG}$ in
  the uniquely ergodic case, see Lemmas~\ref{lem:freq}
  and~\ref{lem:cylinder}. This proves the statement, if
  $\mu(\Xi_{P,U})=0$.  Otherwise, $\nu(P|G \wedge B_{n})$ grows
  unboundedly in $n$.  In particular, $\nu(P|G \wedge B_{n}) \ne 0$
  for $n\ge n_0$.  Thus, we can write
  \begin{equation}
    \label{enlarge}
    \frac{\nu(P^{(\eta)}|G^{(\omega)} \wedge B_{n})}{\vol(B_{n})}
    =
     \frac{\nu(P| G \wedge B_{n})}{\vol(B_{n})} \;
    \frac{\nu(P^{(\eta)}|G^{(\omega)} \wedge B_{n})}{\nu(P|G \wedge B_{n})} 
  \end{equation}
  for $n\ge n_0$. It suffices to show that the second factor on the r.h.s.\ of
  \eqref{enlarge} converges to $\PP_{P}(\eta)$ as $n\to\infty$ for
  $\PP_{G}$-a.a.\ $\omega \in\Omega_{G}$. The latter statement follows from
  the strong law of large numbers, which is applicable due to Kolmogorov's
  criterion \cite{Bau01}. This is obvious, if all translates
  $P_{j}$, $j\in\mathbb{N}$, of $P$ in $G$ are
  pairwise overlap-free and so that colours are assigned in a stochastically
  independent way to different translates.

  Otherwise, we partition the set $\{P_{j}\}_{j\in\mathbb{N}}$ of all
  such translates into a finite number $\Delta$ of (non-empty) subsets
  $\{P_{j}\}_{j\in J_{\alpha}}$, $\alpha =1,\ldots,\Delta$, such that
  there are no mutual overlaps between the translates within each of
  these subsets.  Here we have $\varnothing \neq J_{\alpha}
  \subset\mathbb{N}$ for all $\alpha =1,\ldots,\Delta$,
  $\cup_{\alpha=1}^{\Delta} J_{\alpha} = \mathbb{N}$ and $J_{\alpha}
  \cap J_{\alpha'} = \varnothing$ for all $\alpha \neq \alpha'$.  The
  existence of such a partition may be seen by a graph-colouring
  argument: construct a graph $\mathcal{T}$ such that each translate
  $P_{j}$ defines one point in the vertex set of $\mathcal{T}$.  Two
  vertices are adjacent in $\mathcal{T}$, if the corresponding
  translates overlap.  Clearly, a vertex colouring of $\mathcal{T}$
  (with adjacent vertices having different colours) provides an
  example for the partition that we are seeking. Due to uniform
  discreteness, the degree of any vertex in $\mathcal{T}$ is bounded
  by some number $d_{\mathcal{T}\!,\mathrm{max}} < \infty$. Thus, the
  vertex-colouring theorem \cite{Die05} ensures the existence of such a
  colouring with $\Delta \le 1+ d_{\mathcal{T}\!,\mathrm{max}}$
  different colours. Denoting the number of elements in
  $\{P_{j}\}_{j\in J_{\alpha}}$ with the property $P_{j} \wedge B_{n}
  = P_{j}$ by $\nu_{\alpha}(P |G \wedge B_{n})$, and denoting by
  $\nu_{\alpha}(P^{(\eta)}|G^{(\omega)} \wedge B_{n})$ the analogous
  quantity requiring, in addition, a match of the edge colourings, we
  can write
  \begin{equation}
    \label{decompose}
    \frac{\nu(P^{(\eta)}|G^{(\omega)} \wedge B_{n})}{\nu(P|G \wedge B_{n})} 
    = \sum_{\alpha=1}^{\Delta} \frac{\nu_{\alpha}(P^{(\eta)}|G^{(\omega)}
      \wedge B_{n})}{\nu_{\alpha}(P|G \wedge B_{n})}  \;
    \frac{\nu_{\alpha}(P|G \wedge B_{n})}{\nu(P|G \wedge B_{n})}.
  \end{equation}
  Here we assume $n$ large enough so that $\nu_{\alpha}(P|G \wedge
  B_{n}) >0$ for all $\alpha$. Clearly, those $\alpha$ for which the
  index set $J_{\alpha}$ is finite do not contribute to
  \eqref{decompose} in the macroscopic limit $n\to\infty$, because
  for them the right-most fraction in \eqref{decompose} vanishes in the
  limit. For the remaining $\alpha$'s, the strong law of large numbers
  can be applied and gives
  \begin{equation}
    \lim_{n\to\infty} \frac{\nu_{\alpha}(P^{(\eta)} | G^{(\omega)} \wedge
        B_{n})}{\nu_{\alpha}(P|G \wedge B_{n})} =
      \PP_{P}(\eta)
  \end{equation}
  for $\PP_{G}$-a.a.\ $\omega\in\Omega_{G}$. Therefore we conclude
  \begin{align}
    \lim_{n\to\infty} & \left| \frac{\nu(P^{(\eta)}|G^{(\omega)} \wedge
        B_{n})}{\nu(P|G \wedge B_{n})}  - \PP_{P}(\eta) \right| \nonumber\\ 
    & \le 
    \lim_{n\to\infty} \sum_{\alpha =1}^{\Delta} \left|
      \frac{\nu_{\alpha}(P^{(\eta)} | G^{(\omega)} \wedge 
        B_{n})}{\nu_{\alpha}(P|G \wedge B_{n})}  - \PP_{P}(\eta) \right| \;
    \frac{\nu_{\alpha}(P|G \wedge B_{n})}{\nu(P|G \wedge B_{n})} \nonumber\\
    & \le 
    \lim_{n\to\infty} \sum_{\begin{subarray}{c}\alpha \in\{1,\ldots,\Delta\}: \\
      |J_{\alpha}| =\infty \end{subarray}} \; \left|
      \frac{\nu_{\alpha}(P^{(\eta)} | G^{(\omega)} \wedge 
        B_{n})}{\nu_{\alpha}(P|G \wedge B_{n})}  - \PP_{P}(\eta) \right|
    \nonumber\\     
    & =0.
  \end{align}
  for $\PP_{G}$-a.a.\ $\omega\in\Omega_{G}$, and the lemma follows
  together with \eqref{hof-start}, \eqref{enlarge} and \eqref{freq-conv}. 
\end{proof}

Having established Lemma~\ref{lem:hoflemma}, which is an extension of
\cite[Lemma~3.1]{Hof98} to more general graphs, we can now argue as in
the proof of \cite[Thm.~3.1]{Hof98} to obtain the central result of
this section. 

\begin{theorem}
  \label{theo:perg}
  Let $\cG$ be a set of graphs of finite local complexity and with
  uniformly discrete vertex sets of radius $r>0$. Fix an ergodic
  probability measure $\mu$ on $X_{\cG}$.
  Then there exists a unique ergodic probability measure 
  $\widehat{\mu}$ on $\hXcG$ such that 
  \begin{nummer}
  \item for every pattern $P^{(\eta)}$ of some
  coloured graph in $\hcG$ and every Borel set $U\subseteq\RR^{d}$ with
  diameter $\diam(U) <r$ the relation
    \begin{equation}
      \widehat{\mu}(\widehat{\Xi}_{P^{(\eta)}, U}) = \mu(\Xi_{P, U}) \,
      \PP_{P}(\eta)
    \end{equation}
    holds. 
  \item 
    \label{item-erg}
    for every $\phi\in \mathrm{L}^1(\hXcG,\widehat{\mu})$ the limit
    \begin{equation}
      \label{erg-limit}
      \lim_{n\to\infty} \frac{1}{\vol(B_{n})}  \int_{B_n}\!\d x\;
      \phi(x+G^{(\omega)})  
      =   \int_{\hXcG}
      \!\d\widehat{\mu}(F^{(\sigma)}) \; \phi (F^{(\sigma)})
    \end{equation}
    exist for $\widehat{\mu}$-a.a.\ $G^{(\omega)} \in \hXcG$.
    If $\hXcG$ is even uniquely ergodic and, in
    addition, if $\phi$ is either continuous or a linear combination
    of cylinder functions, then the limit \eqref{erg-limit} exists
    for \emph{all} $G\in X_{\cG}$ and for $\mathbb P_G$-a.a.\
    $\omega\in\Omega_G$.
  \item
    \label{fubini}
    for every $\phi\in \mathrm{L}^1(\hXcG,\widehat{\mu})$ we have
    \begin{equation}
      \int_{\hXcG} \!\d\widehat{\mu}(G^{(\omega)}) \; \phi
      (G^{(\omega)}) =
      \int_{X_{\cG}} \!\d\mu(G) \int_{\Omega_{G}}\!\d\PP_{G}(\omega) \;
      \phi (G^{(\omega)}).
    \end{equation}
  \end{nummer} \qed
\end{theorem}

\begin{remarks} 
\item 
  The corresponding theorem \cite[Thm.~3.1]{Hof98} is a statement
  about Bernoulli site percolation on the Penrose tiling. Our result
  is an extension, which covers both the aperiodic and the periodic
  situations, under weaker assumptions on the base graphs.
\item 
  The asserted uniqueness of the ergodic measure $\widehat{\mu}$ in
  the theorem does not mean that the dynamical system $\hXcG$ is
  uniquely ergodic. It only means that $\widehat{\mu}$ is unique for
  the given ergodic measure $\mu$ on $X_{\cG}$ and the measures
  $\PP_{G}$ on $\Omega_{G}$.
\end{remarks}

%
\section{Finite-range operators on  randomly coloured graphs}
\label{sec:operators}
%

In this section, we consider covariant finite-range operators on
randomly coloured graphs, together with some of their basic ergodic and
spectral properties. We ensure the existence of their integrated
density of states, derive its self-averaging property, and study the
non-randomness of the spectrum. For the spectral-theoretic background,
the reader is referred to \cite{ReSiI,ReSiII}.

For a countable set $\cV$, let $\ell^{2}(\cV)$ be the Hilbert space of
square-summable functions $\psi: \cV \rightarrow\mathbb{C}$ with
canonical scalar product $\langle \cdot, \cdot\rangle$. We denote
the canonical basis in  $\ell^{2}(\cV)$ by $\{\delta_{v}\}_{v\in\cV}$,
that is $\delta_{v}(w) = 1$ if $w=v$ and zero otherwise.

\begin{definition}
  \begin{nummer}
  \item 
    \label{cov-op-def}
    Let $G^{(\omega)}$ be a coloured graph. A bounded linear operator 
    $H_{G^{(\omega)}}$ in $\ell^2(\cV_{G^{(\omega)}})$
    is said to be \emph{covariant of range} $R \in ]0,\infty[$, if  
    \begin{enumerate}
    \item
      $\langle\delta_{x+v}, H_{G^{(\omega)}} \delta_{x+w}\rangle =
      \langle\delta_v, H_{G^{(\omega)}} \delta_w\rangle$ \quad for all $v,w
      \in\cV_{G^{(\omega)}}$ and all $x\in\RR^{d}$ for which $G^{(\omega)}
      \wedge \bigl( B_{R}(x+ v) \cup B_{R}(x+w) \bigr) = G^{(\omega)} \wedge
      \bigl( B_{R}(v) \cup B_{R}(w) \bigr)$,
      \par\smallskip
    \item
      $\langle \delta_v,  H_{G^{(\omega)}}\delta_w\rangle=0$ \quad
      for all $v,w \in \cV_{G^{(\omega)}}$ subject to $|v-w|\ge R$.
    \end{enumerate}
    \smallskip
  \item 
    \label{erg-op-def}
    Let $\cG$ be a set of graphs of finite local complexity and with
    uniformly discrete vertex sets of radius $r$. Fix an ergodic probability
    measure $\widehat{\mu}$ on $\hXcG$. Given any $R\in]0,\infty[$, we call a
    mapping $\widehat{H}: G^{(\omega)} \mapsto H_{G^{(\omega)}}$ from $\hXcG$,
    with values in the set of bounded, self-adjoint operators that are
    covariant of range $R$, a $\widehat{\mu}$-\emph{ergodic self-adjoint
      operator of finite range}.
  \end{nummer}
\end{definition}

\begin{remarks}
\item 
  The covariance condition in the above definition means that
  $H_{G^{(\omega)}}$ is determined on the class of non-equivalent
  $R$-patterns.
\item 
  \label{unibound}
  A $\widehat{\mu}$-ergodic self-adjoint operator of finite range
  $\widehat{H}$ is \emph{uniformly bounded}, in the sense that
  $\sup_{G^{(\omega)} \in \hXcG} \| H_{G^{(\omega)}} \| <
  \infty$, where $\|\cdot\|$ denotes the usual operator norm.  In particular,
  there exists a compact interval $K\subset\RR$, such that for
  $\widehat{\mu}$-almost every $G^{(\omega)} \in \hXcG$ the spectrum of
  $H_{G^{(\omega)}}$ is contained in $K$.
\end{remarks}

The eigenvalue density of $\widehat{H}$ is a quantity of great
interest in applications.

\begin{definition}
  \label{Ndef}
  Let $\widehat{H}$ be a $\widehat{\mu}$-ergodic
  self-adjoint operator of finite range, and fix a
  mollifier as in Corollary~\ref{density}. Then, the \emph{integrated density
    of states} of $\widehat{H}$ is defined as the right-continuous
  distribution function $\RR \rightarrow [0,\varrho]$,
  \begin{equation}
    \label{eq:Ndef}
    E \mapsto N (E) := \int_{\hXcG} \!\d\widehat{\mu}(G^{(\omega)}) \sum_{v\in
      \cV_{G^{(\omega)}}}  \psi(v)\, 
    \langle\delta_v, \Theta(E-H_{G^{(\omega)}}) \delta_v\rangle .
  \end{equation}
  In \eqref{eq:Ndef} we have denoted the right-continuous Heaviside
  unit-step function by $\Theta := \Chi_{[0,\infty[}$ and the
  mollifier $\psi$ was introduced in Corollary~\ref{density}.
\end{definition}

\begin{remarks}
\item 
  \label{integrable}
  The integrand in \eqref{eq:Ndef} is bounded, see
  Remark~\ref{dens-well-def}. Moreover, it is measurable. In fact, the
  map $\hXcG \rightarrow \mathbb{C}$,
  \begin{equation}
    G^{(\omega)} \mapsto f_{\psi}(G^{(\omega)}) := \sum_{v\in
      \cV_{G^{(\omega)}}} \psi(v)\, \langle\delta_v, f(H_{G^{(\omega)}})
    \delta_v\rangle 
  \end{equation}
  is even continuous for all $f \in \mathrm{L}^{\infty}(\RR)$, thanks to the
  continuity of $\psi$.
\item 
  The unique Borel measure $\d N$ on $\RR$ associated with the distribution
  function $N$ has a total mass, given by the vertex density
  $\varrho$, cf.\ Corollary~\ref{density}. Moreover, $\d N$
  is compactly supported, due to the boundedness of $\widehat{H}$.
\item 
  Clearly, the integrated density of states $N$ depends on the choice of
  the ergodic measure $\widehat{\mu}$ on $\hXcG$. However, ergodicity of
  $\widehat{\mu}$ implies that $N$ does not depend on the choice of
  the mollifier $\psi$.  This will become manifest in Theorem~\ref{maclimit}
  below.
\end{remarks}

The integrated density of states of $\widehat{H}$ can also be
characterised in terms of a macroscopic limit.

\begin{theorem}
  \label{maclimit}
  Let $\widehat{H}$ be a $\widehat{\mu}$-ergodic, self-adjoint operator of
  finite range and let $N$ be its integrated density of states
  \eqref{eq:Ndef}, for some choice of the mollifier $\psi$.
  Then there exists a Borel set $\widehat{A} \subseteq\hXcG$ of full
  probability, $\widehat{\mu}(\widehat{A}) =1$, such that
  \begin{equation}
    N(E) = \lim_{n\to\infty} \biggl\{ \frac{1}{\vol(B_{n})} \;
      \sum_{v\in\cV_{G^{(\omega)}} \cap B_{n}} 
     \langle\delta_v, \Theta(E-H_{G^{(\omega)}}) \delta_v\rangle\biggr\}
  \end{equation}
  holds for all $G^{(\omega)} \in \widehat{A}$ and all $E \in\RR$, except for
  the (at most countably many) discontinuity points of $N$. If $\hXcG$
  is uniquely ergodic, then convergence holds even for \emph{all} $G\in
  X_{\cG}$ and $\PP_{G}$-a.a.\ $\omega\in\Omega_{G}$.
\end{theorem}

\begin{proof}[Sketch of the proof]
  The theorem follows from vague convergence of the associated
  measures by standard arguments \cite[Thms.~30.8,
  30.13]{Bau01}. Vague convergence follows in turn from the Ergodic
  Theorem~\ref{item-erg}, and from the identity
  \begin{equation}
    \sum_{v\in\cV_{G^{(\omega)}} \cap B_{n}}  \langle\delta_v,
    f(H_{G^{(\omega)}}) \delta_v\rangle 
    =  \int_{B_n}\!\d x\; f_{\psi}(x+G^{(\omega)})  +
    \mathcal{O}(n^{d-1}),
  \end{equation}
  which is valid for arbitrary $f\in C_{\mathrm{c}}(\RR)$ and arbitrary
  mollifiers $\psi$. The continuous function $f_{\psi}$, associated with $f$,
  was defined in Remark~\ref{integrable}.
\end{proof}

\begin{remark}
  For systems of randomly coloured subgraphs of $\ZZ^{d}$ one can even
  prove uniform convergence in the energy $E$ \cite{LeMu06}. Such a
  result, which is based on ideas in \cite{LeSt05}, has now been
  generalised in \cite{LeVe07}. In particular, it applies also to
  random colourings of graphs with Delone vertex sets. In addition,
  the origin of the discontinuities of $N$ is related to compactly
  supported eigenfunctions and their sizes to equivariant dimensions
  \cite{LeVe07}.

  Statements corresponding to Theorem~\ref{maclimit}
  for systems of uncoloured Delone sets with finite local complexity
  can be found in \cite{Hof95}, \cite[Prop.~4.6]{LeSt03},
  \cite[Thm.~3]{LeSt05} and \cite[Thm.~6.1]{LePe07}. The papers
  \cite{Hof95} and \cite{LeSt05} deal with strictly ergodic systems,
  for which \cite{LeSt05} establish even uniform convergence in the
  energy $E$.
\end{remark}

Next, we relate the set of growth points of the integrated density of
states $N$ to the spectrum of $\widehat{H}$. Given any a self-adjoint
operator $A$, we denote its spectrum by $\spec(A)$ and the essential
part of the spectrum by $\specess(A)$.  We write $\supp (\d N)$ for
the topological support of the probability measure on $\RR$, whose
distribution function is the integrated density of states $N$. The
topological support can be characterised as
\begin{equation}
  \label{top-supp-def}
  \supp (\d N) = \bigg\{ E\in\RR : \int_{]\lambda,\lambda'[} \d N >0
  \text{\quad for all~} \lambda,\lambda'\in \QQ \text{~with~} \lambda <E
  <\lambda'\bigg\} .
\end{equation}

\begin{theorem}
  \label{spec-ids}
  Let $\widehat{H}$ be a $\widehat{\mu}$-ergodic, self-adjoint
  operator of finite range. Then there exists 
  a Borel set $\overline{X} \subseteq X_{\cG}$ with $\mu(\overline{X})=1$ such
  that
  \begin{equation}
    \label{eq:spec-ids}
    \spec(H_{G^{(\omega)}}) =  \specess (H_{G^{(\omega)}}) = \supp (\d
    N)
  \end{equation}
  for all $G\in\overline{X}$ and $\mathbb{P}_{G}$-almost all
  $\omega\in\Omega_{G}$. Moreover, if $X_{G}$ is uniquely ergodic and
  obeys the positive lower frequency condition, see
  Definition~\ref{posfreq}, then the statement holds even with
  $\overline{X} = X_{\cG}$, that is, for \emph{every} $G\in X_{\cG}$.
\end{theorem}

The theorem follows from Lemmas~\ref{N-sub-spec} and~\ref{spec-sub-N}
below.

\begin{remark}
  The most interesting part of Theorem~\ref{spec-ids} is that the statement
  \eqref{eq:spec-ids} holds for all $G\in X_{\cG}$, under the stronger
  assumptions of unique ergodicity and the positive lower frequency
  condition. For uncoloured Delone systems, such a result can be found
  in \cite[Prop.~7.4]{Hof93} and \cite[Lemma~3.6, Thm.~4.3]{LeSt03},
  cf.\ Remarks~\ref{mini-posfreq}.

  The weaker $\mu$-almost sure statement of Theorem~\ref{spec-ids},
  which holds without the additional hypotheses, is analogous to
  results in \cite{LePe07}, where they are proved in the context of
  ergodic groupoids.
\end{remark}

%
\section{Lifshits tails for the Laplacian on bond-percolation graphs}
\label{secLif}

In this section, we study the asymptotics at the lower spectral edge of the
integrated density of states of graph Laplacians, which are associated to
(Bernoulli) bond-percolation subgraphs of a given graph. We will specialise
the general framework of the previous sections in three respects.

First, we set $\cG := \{G_{0}\}$, where $G_{0}$ is some fixed graph of finite
local complexity, with a uniformly discrete vertex set and a bounded degree
sequence. For notational simplicity we write $X$ instead of $X_{\{G_{0}\}}$,
and $\widehat{X}$ instead of $\widehat{X}_{\{G_{0}\}}$. We fix an ergodic
measure $\mu$ on $X$. 

Second, we interpret a Bernoulli bond-percolation subgraph as a
particular randomly coloured graph. To this end, we choose the set of
colours $\AA = \{0,1\}$ and make the identification of a coloured
graph $(G,\omega) \in\widehat{X}$ with the subgraph $(\cV_{G},
\cE_{G}^{(\omega)})$ of $G$ which has the same vertex set as $G$ and
edge set $\cE_{G}^{(\omega)} := \{e\in\cE_{G}: \omega_{e}=1\}$.  We
denote this subgraph again by $G^{(\omega)}$. The probability measure
on $\Omega_{G}$ is given by the $|\cE_{G}|$-fold product measure
$\PP^{p}_{G} = \bigotimes_{e \in \cE_{G}}\PP_{0}^{p}$ of the Bernoulli
measure $\PP_{0}^{p}:=p \boldsymbol{\delta}_{1} +
(1-p)\boldsymbol{\delta}_{0}$ on $\AA$ with parameter $p \in[0,1]$.
Thus, each edge is taken away from $G$ independently of the others
with probability $1-p$. According to Section~\ref{sec:ran-col}, the
ergodic measure $\mu$ on $X$ gives rise to an ergodic measure
$\widehat{\mu}_{p}$ on $\widehat{X}$.

Third, we take the operator $H_{G^{(\omega)}}$ to be the combinatorial or
graph Laplacian $\Delta_{G^{(\omega)}}$ associated with the graph
$G^{(\omega)}$. The graph Laplacian is a covariant operator of range $R=2$ in
the sense of the previous section, as follows easily from its definition in

\begin{lemma}
  Given an arbitrary graph $G=(\cV,\cE)$ with bounded degree sequence, the
  \emph{graph Laplacian} $\Delta_{G}: \ell^{2}(\cV) \rightarrow \ell^{2}(\cV)$,
  \begin{equation}
    (\Delta_{G}\varphi)(v) := \sum_{w \in \cV: \{v,w\} \in \cE}
    [\varphi(v) -\varphi(w)]
  \end{equation}
  for all $\varphi \in\ell^{2}(\cV)$ and all $v \in\cV$, is a bounded,
  self-adjoint and non-negative linear operator. Moreover, zero is an
  eigenvalue of $\Delta_{G}$, if and only if $G$ possesses a finite cluster.
  The multiplicity of the eigenvalue zero is given by the number of finite
  clusters of $G$. \qed
\end{lemma}

We omit the straightforward proof of the lemma and refer to standard
accounts \cite{Chu97,Col98} on spectral graph theory instead. The
mapping $\widehat{X} \ni G^{(\omega)} \mapsto \Delta_{G^{(\omega)}}$
defines a $\widehat{\mu}_{p}$-ergodic, self-adjoint operator of finite
range in the sense of Definition~\ref{erg-op-def}, the \emph{Laplacian
  on bond-percolation graphs associated with} $G_{0}$.

The graph $G_{0}$ need not be connected in what follows. However, it
is crucial that the vertex density $\varrho_{\infty}$ of its infinite
component(s), which was introduced in Corollary~\ref{rho-infty}, is
strictly positive. 

\begin{theorem}
  \label{main}
  Let $G_{0}$ be a graph of finite local complexity, with a uniformly
  discrete vertex set, a bounded degree sequence $d_{\mathrm{max}} :=
  \sup_{v\in\cV} d_{G_{0}}(v) < \infty$, and a maximal edge length 
  $l_\mathrm{max}:=\sup\{|u-v|:\{u,v\}\in \cE\}<\infty$. Assume further 
  that $\varrho_{\infty} > 0$, with respect to some ergodic measure $\mu$ 
  on $X$. Consider the Laplacian on bond-percolation graphs associated
  with $G_{0}$, and let $N$ be its integrated density of states
  \eqref{eq:Ndef}, with respect to the measure $\widehat{\mu}_{p}$.
  If $p \in ]0, \frac{1}{d_{\mathrm{max}}-1}[$, then
  \begin{equation}
    \label{main-eq}
    \lim_{E\downarrow 0}\;\frac{\ln\bigl| \ln [N(E) -N(0)]\bigr|}{\ln E} =
    -1/2 ,
  \end{equation}
  that is, $N$ exhibits a \emph{Lifshits tail} with Lifshits exponent $1/2$ at
  the lower edge of the spectrum.
\end{theorem}

\begin{remark}
  \begin{nummer}
  \item Theorem~\ref{main} follows from Lemmas~\ref{upper}
    and~\ref{lower} below, which provide slightly stronger statements
    as needed to conclude \eqref{main-eq}. The lemmas show that
    Theorem~\ref{main} holds for all edge probabilities $p$, for which
    the cluster-size distribution decays exponentially for
    $\mu$-almost every graph $G\in X$. This means in particular, that
    the validity of Theorem~\ref{main} is limited to the
    \emph{non-percolating phase}, that is $\bigl\{p \in [0,1]:
    \sup_{v\in\cV_{G_{0}}} \PP_{G_{0}}^{p}(|C_{v}| = \infty) =0 \bigr\}$, see
    Section~\ref{sec:percest}.  Exponential decay of the cluster-size
    distribution is guaranteed by the conditions
    $l_\mathrm{max}<\infty$ and $p \in ]0,
    \frac{1}{d_{\mathrm{max}}-1}[$, see Corollary~\ref{together}.
  \item Remark~\ref{rhopos} states sufficient conditions for
    $\varrho_{\infty} >0$.
  \item Theorem~\ref{main} generalises part of Theorem~1.14 in \cite{KiMu06},
    where the special case $G_{0} = \mathbb{L}^{d}$ of the $d$-dimensional
    integer lattice was considered.
  \item The Lifshits exponent $1/2$ in \eqref{main-eq} does not depend on the
    spatial dimension $d$ of the underlying space. This comes from the fact
    that the asymptotics \eqref{main-eq} is determined by the longest linear
    clusters of the percolation graphs $G^{(\omega)} \in \widehat{X}$.
  \end{nummer}
\end{remark}

\begin{lemma}
  \label{upper}
  Let $G_{0}$ be a graph of finite local complexity and with a uniformly
  discrete vertex set. Assume there exists $p_{0}\in ]0,1[$ such that for
  every $p \in [0, p_{0}[$ and $\mu$-almost all $G\in X$ the cluster-size
  distribution for percolation on $G$ decays exponentially, i.e.,
  \begin{equation}
    \label{CSD-decay}
    \mathbb{P}_{G}^{p} \{\omega \in\Omega_{G}: |C_{v}^{(\omega)}| \ge n
    \} \le D(p) \, \exp\{ - \lambda(p) \, n\}
  \end{equation}
  for all $n\in\NN$, where $D(p), \lambda(p) \in ]0,\infty[$ are
  constants that depend on $p$, but are uniform in $G\in X$ and
  $v\in\cV_{G}$. Here $C_{v}^{(\omega)}$ denotes the cluster of the
  graph $G^{(\omega)}$ containing $v\in\mathcal{V}_{G}$. Then
  \begin{equation}
    N(E) - N(0) \le  \varrho D(p)  \exp\{ - \lambda(p)\, E^{-1/2}\}
  \end{equation}
  for all  $p \in [0, p_{0}[$ and all $E>0$. The vertex density $\varrho$ was
  introduced in \eqref{rho-def}.
\end{lemma}

\begin{proof}
  The proof is analogous to that of Lemma~2.7 (Neumann case) in
  \cite{KiMu06}. The block-diagonal form of the Laplacian with respect
  to the cluster structure implies
  \begin{align}
    N(E) -N(0) &= \int_{X}\!\d\mu(G) \, \sum_{v\in\mathcal{V}_{G}}
    \psi(v) \int_{\Omega_{G}} \d\mathbb{P}_{G}^{p}(\omega) \nonumber\\
    & \hspace*{4cm} \times
    \langle\delta_{v}, \bigl[ \Theta (E- \Delta_{C_{v}^{(\omega)}}) -
    \Theta ( \Delta_{C_{v}^{(\omega)}}) \bigr]
    \delta_{v}\rangle
    \nonumber\\
    &  \le\int_{X}\!\d\mu(G) \, \sum_{v\in\mathcal{V}_{G}}
    \psi(v) \int_{\Omega_{G}} \d\mathbb{P}_{G}^{p}(\omega) \,  \Theta
    \bigl(E- E_{1}(C_{v}^{(\omega)})\bigr)  \nonumber\\
    & \hspace*{4cm} \times
    \langle\delta_{v}, \bigl[ \one -
    \Theta ( \Delta_{C_{v}^{(\omega)}}) \bigr]
    \delta_{v}\rangle
    \nonumber\\
    &  \le\int_{X}\!\d\mu(G) \, \sum_{v\in\mathcal{V}_{G}}
    \psi(v) \; \mathbb{P}_{G}^{p} \big\{ \omega\in\Omega_{G} :
     E \ge E_{1}(C_{v}^{(\omega)}) \bigr\} , 
     \label{upper-init}
  \end{align}
  where $E_{1}(C_{v}^{(\omega)})$ denotes the smallest non-zero
  eigenvalue of the Laplacian on the cluster $C_{v}^{(\omega)}$. As a
  particular consequence of the decay \eqref{CSD-decay} of the
  cluster-size distribution, we infer that $|C_{v}^{(\omega)}| <
  \infty$ for all $v\in\cV_{G}$ holds for $\mu$-almost all $G \in X$
  and $\PP_{G}^{p}$-almost all $\omega\in\Omega$.  Hence, Cheeger's
  inequality $E_{1}(C_{v}^{(\omega)}) \ge |C_{v}^{(\omega)}|^{-2}$,
  see e.g.\ Lemma~A.1 in \cite{KhKi06}, can be applied to estimate the
  probability in \eqref{upper-init}, and we get
  \begin{equation}
    \label{upper-final}
    N(E) -N(0) \le  \int_{X}\!\d\mu(G) \, \sum_{v\in\mathcal{V}_{G}}
    \psi(v) \; \mathbb{P}_{G}^{p} \big\{ \omega\in\Omega_{G} :
     |C_{v}^{(\omega)}| \ge E^{-1/2} \bigr\}.     
  \end{equation}
  The lemma now follows from the exponential decay \eqref{CSD-decay} of the
  cluster-size distribution and from \eqref{rho-def}.
\end{proof}

\begin{lemma}
  \label{lower}
  Let $G_{0}$ be a graph of finite local complexity, with
  a uniformly discrete vertex and bounded degree sequence $d_{\mathrm{max}} :=
  \sup_{v\in\cV} d_{G_{0}}(v) < \infty$. Let $p\in ]0,1[$ and $E>0$. Then
  \begin{equation}
    N(E) -N(0) \ge \varrho_{\infty}\, \e^{-2\gamma(p)}  \;
    \frac{\exp\{ - 4 \,\gamma(p) \, E^{-1/2}\}}{2+ 4E^{-1/2}},
  \end{equation}
  where $\gamma(p):= - \ln p - d_{\mathrm{max}} \ln(1-p) >0$.
\end{lemma}

\begin{proof}
  We adapt the strategy of the proof of Lemma~2.9 (Neumann case) in
  \cite{KiMu06}. In the present setting, we have to cope with the additional
  difficulty that vertices in $G$ which are connected can be very far apart in
  the Euclidean metric.  

  Fix $E>0$ arbitrary and let $\{\varepsilon_{j}\}_{j\in\mathbb{N}}$ be a null
  sequence of positive reals, such that $E+\varepsilon_{j}$ is a point of
  continuity of the integrated density of states $N$ for all $j\in\mathbb{N}$.
  Then the right-continuity of $N$, the Ergodic Theorem~\ref{maclimit} and the
  isotony of $N$ imply
  \begin{multline}
    \label{low-start}
    N(E) -N(0)  = \lim_{j\to\infty} [ N(E+\varepsilon_{j})  -N(0) ] \\
    \ge \limsup_{n\to\infty} \frac{1}{\vol(B_{n})} \;
  \sum_{v\in\cV_{G}} \Chi_{B_{n}}(v) \;\langle\delta_{v}, [
  \Theta (E- \Delta_{G^{(\omega)}}) - \Theta (\Delta_{G^{(\omega)}})]
  \delta_{v}\rangle 
  \end{multline}
  for $\widehat{\mu}_{p}$-almost all graphs $G^{(\omega)} \in\widehat{X}$.
  Since $\Delta_{G^{(\omega)}}$ is a direct sum of the Laplacians of the
  clusters of $G^{(\omega)}$, and since this is also true for functions of the
  Laplacian, it follows that the trace on the right-hand side of
  \eqref{low-start} can be bounded from below by throwing away all
  contributions from branched clusters in that sum,
  \begin{multline}
    \label{throw-away}
      \sum_{v\in\cV_{G}} \Chi_{B_{n}}(v) \;\langle\delta_{v}, [
    \Theta (E- \Delta_{G^{(\omega)}}) - \Theta (\Delta_{G^{(\omega)}})]
     \delta_{v}\rangle  \\
    \ge \sum_{l=2}^{\infty}
    Z_{B_{n}}^{G^{(\omega)}} (\mathcal{L}_{l}) \,
    \bigl\langle \delta_{1}, [  \Theta (E- \Delta_{\mathcal{L}_{l}}) - \Theta
    (\Delta_{\mathcal{L}_{l}}) ] \delta_{1}\bigr\rangle.
  \end{multline}
  Here $\mathcal{L}_{l}$ denotes a linear chain (i.e.\ non-branched
  and cycle-free cluster) with $l$ vertices, $Z_{B_{n}}^{G^{(\omega)}}
  (\mathcal{L}_{l})$ the number of such chains in the percolation
  graph $G^{(\omega)}$, subject to the condition that at least one of
  its end-vertices lies in the ball $B_{n}$. The symbol $\delta_{1}$
  denotes the canonical basis vector in $\ell^{2}(\{1,\ldots,l\})$
  corresponding to one end-vertex of $\mathcal{L}_{l}$ (by symmetry
  reasons it does not matter which end-vertex).

  The spectral representation of $\Delta_{\mathcal{L}_{l}}$ with $l \ge 2$ is
  explicitly known, for example, by mapping the problem to that of a cycle
  graph with $2l$ vertices. The eigenvalues turn out to be
  \begin{equation}
    E_{k}(\mathcal{L}_{l}) := 4 \bigl(\sin(\pi k/2l)\bigr)^{2} , \qquad\quad
    k=0,\ldots l,
  \end{equation}
  and the components of the corresponding normalised eigenvectors
  $\varphi_{k}$ in the canonical basis are given by
  $\langle\delta_{j},\varphi_{0}\rangle := l^{-1/2}$ and
  \begin{equation}
    \langle\delta_{j}, \varphi_{k}\rangle :=  (
    2/l)^{1/2} \cos \bigl( \pi
    \tfrac{k}{l}\,(j- \tfrac{1}{2})\bigr), \qquad\quad
    k=1,\ldots l,
  \end{equation}
  where $j=1,\ldots, l$.  We observe that
  \begin{equation}
    \label{op-bound}
    \Theta (E- \Delta_{\mathcal{L}_{l}}) - \Theta
    (\Delta_{\mathcal{L}_{l}}) \ge \Theta (E- E_{1}(\mathcal{L}_{l})) \;
    \varphi_{1}  \otimes \varphi_{1},
  \end{equation} 
  where the dyadic product is the projector on the eigenspace generated by
  $\varphi_{1}$, and that $E_{1}(\mathcal{L}_{l}) \le 10/l^{2}$ and
  $|\langle \delta_{1}, \varphi_{1}\rangle|^{2} \ge l^{-1}$.  Therefore we
  obtain
  \begin{align}
    \label{number-chain}
    N(E) -N(0)  & \ge \limsup_{n\to\infty} \sum_{l=2}^{\infty}
    \Theta (E- 10/l^{2}) \;\frac{1}{l} \;
    \frac{Z_{B_{n}}^{G^{(\omega)}} (\mathcal{L}_{l})}{\vol(B_{n})} \nonumber\\
    & \ge \frac{1}{l(E)} \; \limsup_{n\to\infty}
    \frac{Z_{B_{n}}^{G^{(\omega)}} 
      (\mathcal{L}_{l(E)})}{\vol(B_{n})}
  \end{align}
  for $\widehat{\mu}_{p}$-almost all graphs $G^{(\omega)}\in\widehat{X}$ with
  $l(E):= \inf\{l\in\mathbb{N}\setminus\{1\}: E -10/l^{2} \ge 0\}$.

  The quantity 
  \begin{equation}
    g_{v}(G^{(\omega)}) := \left\{
      \begin{array}{l@{\quad}l}
        1, & \mbox{if a vertex $v\in\mathcal{V}_{G}$ is an end-vertex
          of a linear chain}\\[-.2ex]
        & \mbox{with $l(E)$ vertices in $G^{(\omega)}$}\\[.5ex]
        0, & \mbox{otherwise}
      \end{array}\right.
  \end{equation}
  helps to rewrite the right-hand side of \eqref{number-chain}, so that
  \begin{align}
    N(E) -N(0) & \ge \frac{1}{2l(E)} \limsup_{n\to\infty} \frac{1}{\vol(B_{n})}\;
    \sum_{v\in\mathcal{V}_{G}} \Chi_{B_{n}(v)} \, g_{v}(G^{(\omega)})
    \nonumber\\
    & = \frac{1}{2l(E)} \limsup_{n\to\infty} \frac{1}{\vol(B_{n})}
    \int_{B_{n}} \!\d x\, \sum_{v\in\mathcal{V}_{x+G}} \psi(v)
    g_{v}(x+G^{(\omega)})
  \end{align}
  for $\widehat{\mu}_{p}$-almost all graphs $G^{(\omega)} \in \widehat{X}$.
  The Ergodic Theorem~\ref{item-erg},~\itemref{fubini} now implies 
  \begin{equation}
    N(E) -N(0) \ge  \frac{1}{2l(E)} \int_{X}\!\d\mu(G)
    \sum_{v\in\mathcal{V}_{G}} \psi(v)
    \int_{\Omega_{G}} \! \d\mathbb{P}^{p}_{G}(\omega) \, g_{v}(G^{(\omega)}).
  \end{equation}
  We recall that $\Chi_{\cV_{{}G,\infty}}(v) =1$, if $v$ belongs to an infinite
  component of $G$, and zero otherwise. Then we have for all $G \in X$ the
  crude elementary combinatorial estimate
  \begin{equation}
    \label{chain-exist}
    \int_{\Omega_{G}}\d\mathbb{P}^{p}_{G}(\omega)\,  g_{v}(G^{(\omega)})
    \ge 2 p^{l(E)} (1-p)^{l(E)d_{\mathrm{max}}} \,\Chi_{\cV_{G,\infty}}(v)
  \end{equation}
  for the probability that a given vertex appears as the end vertex of a
  linear chain with $l(E)$ vertices. Here we have also used Lemma~\ref{dmax}.
  The estimate \eqref{chain-exist} yields
  \begin{equation}
    N(E) -N(0) \ge \frac{\e^{-l(E) \gamma(p)}}{l(E)} \int_{X}\!\d\mu(G)
    \sum_{v\in\mathcal{V}_{G}} \psi(v) \, \Chi_{\cV_{G,\infty}}(v) 
    = \varrho_{\infty} \; \frac{\e^{-l(E) \gamma(p)}}{l(E)} .
  \end{equation}
  Making use of $l(E) < 2 + 4 E^{-1/2}$, we obtain the assertion of the lemma.
\end{proof}

%
\section{Proof of Theorem~\ref{spec-ids}} 
\label{sec:spectrum}
%

Lemmas~\ref{N-sub-spec} and~\ref{spec-sub-N} in this section provide
the proof of Theorem~\ref{spec-ids}.  We begin with a standard result
on the non-randomness of the spectrum.

\begin{lemma}
  Let $\widehat{H}$ be a $\widehat{\mu}$-ergodic, self-adjoint
  operator. Then there exists a Borel set $\Sigma \subseteq \RR$, such
  that $\spec (H_{G^{(\omega)}}) =\Sigma$ for $\widehat{\mu}$-almost
  all $G^{(\omega)} \in\hXcG$. \qed
\end{lemma}

\begin{remarks}
\item 
  \label{proj-det}
  The proof of the lemma is classical in the theory of random
  Schr\"odinger operators, see e.g.\ \cite{CaLa90,PaFi92}. The central
  point in the argument is that for every Borel set $I\subseteq\RR$
  the measurable function
  \begin{equation}
    \label{tdef}
    \hXcG \ni G^{(\omega)} \mapsto t_{I}(G^{(\omega)}) :=
  \tr \Chi_{I}(H_{G^{(\omega)}})
  \end{equation}
  is invariant under translations of
  $G^{(\omega)}$ by arbitrary $a\in\RR^{d}$, and hence
  $\widehat{\mu}$-almost surely constant by ergodicity. 
\item 
  Standard arguments extend the non-randomness also to the
  Lebesgue components of the spectrum \cite{CaLa90,PaFi92}. These
  arguments are taken up in \cite{LePe07} in the context of ergodic
  groupoids.
\end{remarks}

\begin{lemma}
  \label{lemma-equiv}
  Let $\widehat{H}$ be a $\widehat{\mu}$-ergodic, self-adjoint
  operator of finite range, and let $N$ be its integrated density of
  states.  For a given open interval $I:= ]\lambda ,\lambda'[$,
  where $\lambda,\lambda' \in\RR$ with $\lambda <\lambda'$, consider the
  following statements.
  \begin{enumerate}
  \item[(i)] 
    ~~$\displaystyle \int_{I} \d N > 0$.
  \item[(ii)] 
    there exists a Borel set $\overline{X} \subseteq X_{\cG}$ with
    $\mu(\overline{X}) = 1$ \textup{(}$\overline{X} = X_{\cG}$, if
    $X_{\cG}$ is uniquely 
    ergodic\textup{)} such that $t_{I}(G^{(\omega)}) = \infty$ for
    all $G\in \overline{X}$ and for $\PP_{G}$-almost all
    $\omega\in\Omega_{G}$.
  \item[(iii)] 
    there exists a Borel set $\overline{X} \subseteq X_{\cG}$ with
    $\mu(\overline{X}) = 1$  \textup{(}$\overline{X} = X_{\cG}$, if
    $X_{\cG}$ is uniquely 
    ergodic\textup{)} such that $t_{I}(G^{(\omega)}) > 0$ for
    all $G\in \overline{X}$ and for $\PP_{G}$-almost all
    $\omega\in\Omega_{G}$.
  \item[(iv)]
    there exists a Borel set $\underline{X} \subseteq X_{\cG}$ with
    $\mu(\underline{X}) >0$ such that $t_{I}(G^{(\omega)}) > 0$
     for all $G\in \underline{X}$ and all $\omega$ in some subset of
     $\Omega_{G}$ that has a positive $\PP_{G}$-measure.
   \item[(v)]
     there exists $G\in X_{\cG}$ and $\omega\in\Omega_{G}$ such that
     $t_{I}(G^{(\omega)}) > 0$. 
  \end{enumerate}
  Then the implications
  \begin{equation}
    \mathrm{(i)}   \; \Longleftrightarrow \; 
    \mathrm{(ii)}  \; \Longleftrightarrow \;  
    \mathrm{(iii)} \; \Longleftrightarrow \;
    \mathrm{(iv)}  \; \Longrightarrow \; \mathrm{(v)}
  \end{equation}
  hold.  Moreover, if $X_{G}$ is uniquely ergodic and
  obeys the positive lower frequency condition, then 
  \begin{equation}
     \mathrm{(v)}  \; \Longrightarrow \; \mathrm{(i)}
  \end{equation}
  holds, too.
\end{lemma}

\begin{proof}
  (i) $\Rightarrow$ (ii): \quad Using the Ergodic
  Theorem~\ref{maclimit}, we deduce the existence of a Borel set
  $\overline{X} \subseteq X_{\cG}$ with $\mu(\overline{X}) =1$ (and $\overline{X} =
  X_{\cG}$, if $X_{\cG}$ is uniquely ergodic) such that 
  \begin{align}
    0 < \int_{I}\!\d N &= \int_{X_{\cG}}\!\d\mu(\widetilde{G})
    \int_{\Omega_{\widetilde{G}}} \!\d\PP_{\widetilde{G}}(\widetilde{\omega})  \;
    \sum_{v\in\cV_{\widetilde{G}}} \psi(v) \, \langle\delta_{v},
    \Chi_{I}(H_{\widetilde{G}^{(\widetilde{\omega})}}) \delta_{v} \rangle \nonumber\\
    &= \lim_{n\to\infty} \frac{1}{\vol(B_{n})} \; \sum_{v\in \cV_{G}
    \cap B_{n}}  \langle\delta_{v},
    \Chi_{I}(H_{G^{(\omega)}}) \delta_{v} \rangle
  \end{align}
  for all $G\in \overline{X}$ and for $\PP_{G}$-almost all
  $\omega\in\Omega_{G}$. Hence,
  \begin{equation}
    t_{I}(G^{(\omega)}) = \lim_{n\to\infty} \sum_{v\in \cV_{G}
    \cap B_{n}}  \langle\delta_{v},
    \Chi_{I}(H_{G^{(\omega)}}) \delta_{v} \rangle =
    \infty
  \end{equation}
  for those $G^{(\omega)}$.

  (ii) $\Rightarrow$ (iii) $\Rightarrow$ (iv) $\Rightarrow$ (v): \quad
  These implications are obvious.

  (iv) $\Rightarrow$ (i): \quad Assume (i) is wrong. Then we have
  \begin{equation}
    \label{zerostart}
    0 =   \int_{I}\!\d N
    = \int_{X_{\cG}}\!\d\mu(G) \int_{\Omega_{G}}\!\d\PP_{G}(\omega)
    \sum_{v\in\cV_{G}} \psi(v) \, \langle\delta_{v},
    \Chi_{I}(H_{G^{(\omega)}}) \delta_{v}\rangle.
  \end{equation}
  Note that, by translation invariance of the measure $\mu$, we can
  replace the mollifier $\psi$ in \eqref{zerostart} by $\psi_{a}:=
  \psi(\cdot -a)$, for any translation vector $a\in\RR^{d}$.
  Consequently, there exists a Borel set $\overline{X}_{a} \subseteq
  X_{\cG}$ with $\mu(\overline{X}_{a}) = 1$, such that for all
  $G\in\overline{X}_{a}$ there is $\overline{\Omega}_{G,a} \subseteq \Omega_{G}$
  measurable with $\PP_{G}(\overline{\Omega}_{G,a}) = 1$, such that for all
  $\omega\in \overline{\Omega}_{G,a}$ we have $ \langle\delta_{v},
  \Chi_{I} (H_{G^{(\omega)}}) \delta_{v}\rangle =0$,
  for all $v\in\cV_{G} \cap \supp(\psi_{a})$. 
  We can choose a countable set $M$ of translation vectors, such that
  for every $v\in\cV_{G}$ there exists $a\in M$ with $v\in \supp\psi_{a}$. 
  Next, we define $\overline{X} := \bigcap_{a\in M} \overline{X}_{a}$ and
  $\overline{\Omega}_{G} := \bigcap_{a\in M} \overline{\Omega}_{G,a}$ for all
  $G\in \overline{X}$ so that  $\mu(\overline{X}) = 1$ and
  $\PP_{G}(\overline{\Omega}_{G}) = 1$ for all $G\in\overline{X}$. Now, we get
  $\langle\delta_{v}, \Chi_{I} (H_{G^{(\omega)}})
  \delta_{v}\rangle =0$ for all $G\in\overline{X}$, all $\omega\in
  \overline{\Omega}_{G}$ and all $v\in\cV_{G}$. In other words,
  $\Chi_{I}(H_{G^{(\omega)}}) = 0$ for all
  $G\in\overline{X}$ and all $\omega\in \overline{\Omega}_{G}$. This contradicts
  (iv) so that the implication is proven. 

  (v)~$\Rightarrow$~(i): \quad By assumption, $I$ contains a spectral
  value of $H_{G^{(\omega)}}$. Therefore there exists $E\in I$,
  $\delta \in ]0, \varepsilon[$, where $\varepsilon := \dist(E, \RR
  \setminus I)$, and $\varphi\in\ell^{2}(\cV_{G})$, $\varphi \neq 0$,
  such that
  \begin{equation}
    \label{spec-cond}
    \| (H_{G^{(\omega)}} -E) \varphi\| < \delta \|\varphi\|.
  \end{equation}
  Since $H_{G^{(\omega)}}$ is bounded, we can even assume that
  $\varphi$ has compact support. Furthermore, since $\widehat{H}$ is
  of finite range $R$, we choose a compact subset $K_{0}$ of $\RR^{d}$
  such that $\dist\bigl(\supp(\varphi), \RR^{d} \setminus K_{0}\bigr)
  > 2R$.  We write $P_{0} := G \wedge K_{0}$ for the corresponding
  pattern of $G$. From the positive lower frequency condition we infer that
  the copies $P_{j} := a_{j} +P_{0}$, $j\in\NN$, $a_{j} \in\RR^{d}$,
  of this pattern in $G$ occur with a positive lower frequency. For any
  given $P_{j}$ there is a maximum number (which is uniform in
  $j$) of other copies
  $P_{j'}$ with which $P_{j}$ can overlap. Therefore, from now on, we
  will pass to a subsequence of the sequence
  $\{P_{j}\}_{j\in\NN_{0}}$ such
  that none of the patches in the subsequence overlap, and still
  \begin{equation}
    \label{posfreqcond}
    \liminf_{n\to\infty} \frac{\widetilde{\nu}(P_{0}|G\wedge
      B_{n})}{\vol(B_{n})} =:\gamma >0. 
  \end{equation}
  Here the symbol $\tilde{\nu}$ is used instead of $\nu$ to indicate
  that it is only the patterns in the subsequence which are counted 
  in $G\wedge B_{n}$. We denote the subsequence again by
  $\{P_{j}\}_{j\in\NN_{0}}$ and introduce the translated functions
  $\varphi_{j} := \varphi( \cdot -a_{j})$ for $j\in\NN$. They form an
  orthogonal sequence, because $\supp(\varphi_{j}) \subset K_{j} :=
  a_{j} + K_{0}$, and the $K_{j}$'s are pairwise disjoint.  Next we
  have to ensure that the colouring of $G^{(\omega)}$ in $K_{0}$ is
  also repeated in sufficiently many of the $K_{j}$. The events
  \begin{equation}
    A_{j} := \bigl\{ \overline{\omega} \in \Omega_{G} :
    \overline{\omega}_{a_{j} +v} = \omega_{v} \text{~for all~} v\in
    \cV_{G} \cap K_{0}\bigr\},
  \end{equation}
  $j \in\NN_{0}$, where $a_{0} :=0$, are all independent and
  $\PP_{G}(A_{j}) = \PP_{G}(A_{0}) >0$ for all $j\in\NN$. Thus, the
  strong law of large numbers gives
  \begin{equation}
    \label{slln}
    \lim_{n\to\infty} \, \frac{1}{\widetilde{\nu}(P_{0}|G\wedge B_{n})}
    \sum_{j\in\NN_{0} : 
      K_{j} \subset B_{n}} \Chi_{A_{j}}(\overline{\omega}) =
    \PP_{G}(A_{0}) >0
  \end{equation}
  for $\PP_{G}$-almost all $\overline{\omega} \in\Omega_{G}$. 

  We conclude from the Ergodic Theorem~\ref{maclimit} for uniquely
  ergodic systems that
  \begin{equation}
    \label{erg-start}
    \int_{I}\!\d N = \lim_{n\to\infty} \, \frac{1}{\vol(B_{n})}
    \sum_{v\in \cV_{G} \cap B_{n}} \langle\delta_{v},
    \Chi_{I}(H_{G^{(\overline{\omega})}}) \delta_{v} \rangle,
  \end{equation}
  for the given $G\in X_{\cG}$ (for which (v) holds), and for all
  $\overline{\omega}$ in some measurable set $\overline{\Omega}_{G} \subset
  \Omega_{G}$ of full $\PP_{G}$-measure. We choose $\overline{\Omega}_{G}$
  such that \eqref{slln} holds for all
  $\overline{\omega}\in\overline{\Omega}_{G}$, too. 
  Then, we rewrite
  \begin{align}
     \sum_{v\in \cV_{G} \cap B_{n}} \langle\delta_{v},
    \Chi_{I}(H_{G^{(\overline{\omega})}}) \delta_{v} \rangle
    & = \tr_{\ell^{2}(\cV_{G})} \Big\{ \Chi_{B_{n}}
    \Chi_{I}(H_{G^{(\overline{\omega})}}) \Chi_{B_{n}} \Big\} \nonumber\\
    & \ge  \sum_{j\in\NN_{0} : K_{j} \subset B_{n}}
    \frac{1}{\|\varphi_{j}\|^{2}} \,   \langle\varphi_{j},
    \Chi_{I}(H_{G^{(\overline{\omega})}})
    \varphi_{j}\rangle   \nonumber\\
    & \ge  \sum_{j\in\NN_{0} : K_{j} \subset B_{n}}
      \frac{\Chi_{A_{j}(\overline{\omega})}}{\|\varphi_{j}\|^{2}} \,
      \langle\varphi_{j}, \Chi_{I}(H_{G^{(\overline{\omega})}})
    \varphi_{j}\rangle
  \end{align}
  and observe
  \begin{align}
    \langle\varphi_{j}, \Chi_{I}(H_{G^{(\overline{\omega})}})
    \varphi_{j}\rangle  & = \|\varphi_{j}\|^{2} - 
    \| \Chi_{\RR\setminus I}(H_{G^{(\overline{\omega})}}) \varphi_{j}
    \|^{2} \nonumber\\
    &\ge \|\varphi_{j}\|^{2} - \varepsilon^{-2} \int_{\RR\setminus I}
    \! \d\zeta_{G^{(\overline{\omega})},j}(E')\; (E'-E)^{2} \nonumber\\
    &\ge \|\varphi_{j}\|^{2} - \varepsilon^{-2} \|
    (H_{G^{(\overline{\omega})}} -E)\varphi_{j}\|^{2},
  \end{align}
  where $\zeta_{G^{(\overline{\omega})},j} := \langle\varphi_{j},
  \Chi_{\bullet}(H_{G^{(\overline{\omega})}}) \varphi_{j}\rangle$ is the
  projection-valued spectral measure for $H_{G^{(\overline{\omega})}}$ and
  the vector $\varphi_{j}$.
  For $\overline\omega \in A_{j}$, we have $-a_{j} + (G^{\overline{\omega}}
  \wedge K_{j}) = G^{{\omega}} \wedge K_{0}$. Thus, covariance of
  $\widehat{H}$ and \eqref{spec-cond} imply
  \begin{equation}
    \label{covariance}
    \| (H_{G^{(\overline{\omega})}} -E)\varphi_{j}\| = 
    \| (H_{G^{(\omega)}} -E)\varphi\| < \delta \|\varphi\| =
    \delta\|\varphi_{j}\| .
  \end{equation}
  Combining \eqref{erg-start} -- \eqref{covariance}, we get
  \begin{align}
       \int_{I}\!\d N & \ge 
       (1-\delta^{2}/\varepsilon^{2}) \liminf_{n\to\infty} \bigg\{
       \frac{\widetilde{\nu}(P_{0}|G\wedge B_{n})}{\vol(B_{n})}
       \sum_{j\in\NN_{0} : K_{j} \subset B_{n}}
       \frac{\Chi_{A_{j}}(\overline{\omega})}{\widetilde{\nu}(P_{0}|G\wedge
         B_{n}) } \biggr\}\nonumber\\
       & \ge (1-\delta^{2}/\varepsilon^{2}) \, \gamma \, \PP_{G} (A_{0}) >0.
  \end{align}
  The last inequality uses the positive lower frequency condition
  \eqref{posfreqcond} and the strong law of large numbers
  \eqref{slln}. This completes the proof.
\end{proof}

The next two lemmas provide the proof of Theorem~\ref{spec-ids}. They
are consequences of the previous lemma.

\begin{lemma}
  \label{N-sub-spec}
  Let $\widehat{H}$ be a $\widehat{\mu}$-ergodic, self-adjoint
  operator of finite range and let $N$ be its integrated density of
  states. Then there exists a Borel set $\overline{X} \subseteq X_{\cG}$
  with $\mu(\overline{X}) =1$ such that
  \begin{equation}
    \supp(\d N) \subseteq \specess(H_{G^{(\omega)}})
  \end{equation}
  for all $G\in\overline{X}$ and $\PP_{G}$-almost all $\omega\in\Omega_{G}$.
  If $X_{\cG}$ is even uniquely ergodic, then the statement holds with
  $\overline{X} = X_{\cG}$.
\end{lemma}

\begin{proof}
  By taking countable intersections of sets of full measure, we infer
  from the implication  (i)~$\Rightarrow$~(ii) in
  Lemma~\ref{lemma-equiv}, that there exists a Borel set $\overline{X}
  \subseteq X_{\cG}$ with $\mu(\overline{X}) =1$ ($\overline{X} = X_{\cG}$, if
  $X_{\cG}$ is uniquely ergodic), such that for all
  $G \in \overline{X}$ there exists $\overline{\Omega}_{G} \subset \Omega_{G}$
  measurable with $\PP_{G}(\overline{\Omega}_{G}) =1$, such that for all
  $\omega\in\overline{\Omega}_{G}$ and all  
  $\lambda,\lambda' \in\QQ$ with $\lambda <\lambda'$ we have 
  \begin{equation}
    \int_{]\lambda,\lambda'[}\!\d N > 0 \quad\Longrightarrow\quad
    t_{]\lambda,\lambda'[} (G^{(\omega)}) =\infty.
  \end{equation}
  Thus, we conclude from the characterisation \eqref{top-supp-def}
  that
  \begin{equation}
    \supp (\d N) \subseteq \big\{  E \in\RR:
    t_{]\lambda,\lambda'[}(G^{(\omega)}) =\infty \text{~for all~}
    \lambda,\lambda'\in\QQ \text{~with~} \lambda <E<\lambda'\big\}
  \end{equation}
  for  all  $G\in\overline{X}$ and all $\omega\in\overline{\Omega}_{G}$. But
  this is the claim.
\end{proof}

\begin{lemma}
  \label{spec-sub-N}
  Let $\widehat{H}$ be a $\widehat{\mu}$-ergodic, self-adjoint
  operator of finite range, and let $N$ be its integrated density of
  states. Then there exists a Borel set $\overline{X} \subseteq X_{\cG}$
  with $\mu(\overline{X}) =1$, such that
  \begin{equation}
    \label{eq:specN}
    \spec(H_{G^{(\omega)}}) \subseteq \supp(\d N)
  \end{equation}
  for all $G\in\overline{X}$ and $\PP_{G}$-almost all
  $\omega\in\Omega_{G}$.
  If $X_{\cG}$ is even uniquely ergodic and if the positive lower frequency
  condition holds, then the statement holds with
  $\overline{X} = X_{\cG}$ and for all $\omega\in\Omega_{G}$.  
\end{lemma}

\begin{proof}
  Recall that
  \begin{equation}
    \label{spec-rep}
    \spec(H_{G^{(\omega)}}) =  \big\{  E \in\RR:
    t_{]\lambda,\lambda'[}(G^{(\omega)}) >0 \text{~for all~}
    \lambda,\lambda'\in\QQ \text{~with~} \lambda <E<\lambda'\big\}.
  \end{equation}
  For uniquely ergodic systems that obey the positive lower frequency
  condition, we deduce the desired inclusion \eqref{eq:specN} -- valid
  for all $G\in X_{\cG}$ and all $\omega\in\Omega_{G}$ -- directly from
  the implication (v)~$\Rightarrow$~(i) in Lemma~\ref{lemma-equiv}.

  In the general case, we argue that there exists a Borel set
  $\widehat{A} \subseteq \hXcG$ with $\widehat{\mu}(\widehat{A}) =1$,
  such that for all $\lambda,\lambda' \in\QQ$ with $\lambda <\lambda'$
  the implication 
  \begin{equation*}
    t_{]\lambda,\lambda'[}(G^{(\omega)}) > 0 \text{~~for some~}
    G^{(\omega)} \in \widehat{A} \quad \Longrightarrow \quad
     t_{]\lambda,\lambda'[}(G^{(\omega)}) > 0 \text{~for all~}
    G^{(\omega)} \in \widehat{A}
   \end{equation*}
  holds, confer Remark~\ref{proj-det}.
  Hence, for $G^{(\omega)} \in\widehat{A}$, we deduce the assertion from
  \eqref{spec-rep} and the implication
  (ii)~$\Rightarrow$~(i) in Lemma~\ref{lemma-equiv}.
\end{proof}

%
\section{Percolation estimates} 
\label{sec:percest}
%

In this final section, we establish the percolation estimates that
guarantee exponential decay of the cluster-size distribution on rather
general graphs. Corollary~\ref{together} provides a rough criterion
for this. It is used in the proof of Theorem~\ref{main} to ensure the
applicability of Lemma~\ref{upper}. For this reason, the results of
this section must be valid without any assumptions on the automorphism
group of the graph. So far, this has prevented us from extending the
results of this section to higher values of the percolation
probability up to the critical value. This is a challenging open
problem, see also the discussion in \cite{Hof98}. 
Assuming quasi-transitivity, stronger results were obtained recently in
\cite{AnVe07b}.

\smallskip

First we give a simple, but crude lower bound for the critical
probability of Bernoulli bond percolation on an infinite connected
graph $G=(\cV,\cE)$ with bounded degree sequence.  Let
$\theta_{G,v}(p) := \mathbb P_G^p(|C_v|=\infty)$ denote the
probability that the open cluster containing $v\in\cV$ is infinite.
The critical probability $p_c(G) := \sup\{p \in
[0,1]:\theta_{G,v}(p)=0\}$ is independent of $v\in\cV$, as follows
from the FKG inequality, see \cite[Thm.~2.8]{Gri99} or \cite{Hof98}
(recall that we consider only graphs with countable vertex sets). By
standard reasoning \cite[Thm~1.10]{Gri99}, we have the following
elementary lower estimate for $p_c(G)$.

\begin{lemma}
  \label{combi}
  Let $G=(\cV,\cE)$ be an infinite connected graph.
  If $G$ has maximal vertex degree $d_\mathrm{max} := \sup_{v\in\cV} d_{G}(v)
  \in \mathbb{N} \setminus\{1\}$, then
  \begin{equation}
    p_c(G)\ge\frac{1}{ d_\mathrm{max}-1}.
  \end{equation}
\end{lemma}

\begin{proof}
  Let $\sigma_{v}(n)$ denote the number of $n$-step self-avoiding walks on
  $G$, starting from $v\in V$. Since a self-avoiding walk must not return to
  its previous position when performing a single step, we obtain
  $\sigma_{v}(n)\le d_\mathrm{max} ( d_\mathrm{max}-1)^{n-1}$. Let
  $W_{v}^{(\omega)}(n)$ denote the number of such walks in the percolation
  subgraph $G^{(\omega)}$ of $G$. Since every such walk is open with
  probability $p^n$ in any percolation subgraph, we get for 
  its expectation $\int_{\Omega_{G}}\mathbb{P}_{G}^{p}(\d\omega)\,
  W_{v}^{(\omega)}(n) = p^n \sigma_{v}(n)$.  Note that, if the vertex $v$ belongs
  to an infinite cluster, there are open self-avoiding walks of arbitrary
  length emanating from $v$. Thus, we have for all $n\in\mathbb N$ the
  estimate
  \begin{multline}
    \theta_{G,v}(p)\le \mathbb P_G^p\bigl\{ \omega\in\Omega_{G}:
    W_{v}^{(\omega)}(n)\ge 1 \bigr\} 
    \le \int_{\Omega_{G}}\!\mathbb{P}_{G}^{p}(\d\omega)\,
  W_{v}^{(\omega)}(n) \\
    =  p^n\sigma_{v}(n) \le\frac{ d_\mathrm{max}}{ d_\mathrm{max}-1}
    \;[p(d_\mathrm{max}-1)]^n.  
  \end{multline}
  This implies $\theta_{G,v}(p)=0$, if $p( d_\mathrm{max}-1)<1$, and we obtain
  the assertion of the lemma.
\end{proof}

The behaviour in the subcritical phase $p<p_c(G)$ can be inferred 
from asymptotic properties of the event
\begin{equation}
  A_{v}(n) := \bigl\{\omega \in\Omega_{G}:  v
  \stackrel{G^{(\omega)}}{\longleftrightarrow} B_{n}(v)^c \bigr\}
\end{equation}
that there exists a path from $v\in \mathcal V$ to the complement of the ball
of radius $n$ around $v$. For fixed $p$, denote by $g_{G,v}^{p}(n) := \mathbb
P_{G}^{p}\bigl( A_{v}(n)\bigr)$ the probability of the event $A_{v}(n)$.
The following lemma states that $g_{G,v}^{p}(n)$ decays exponentially in $n$,
if the percolation probability $p$ is small enough. Its proof is analogous to
that of the previous lemma.

\begin{lemma}
  \label{expsmallp}
  Let a graph $G=(\cV,\cE)$ be given. Assume that $G$ has maximal
  vertex degree $d_\mathrm{max} \in \mathbb{N} \setminus\{1\}$ and
  maximal edge length $l_\mathrm{max}:=\sup\{|u-v|:\{u,v\}\in
  \cE\}<\infty$. Then, for every $p\in ]0,1]$ there exists a real
  number $\psi(p)$, such that the probability $g_{G,v}^{p}(n)$ of the
  event $A_{v}(n)$ satisfies
  \begin{equation}
    g_{G,v}^{p}(n)\le 2\, \e^{-n\psi(p)} 
  \end{equation}
  for all $n\in\mathbb N$, uniformly in $v\in\cV$. A possible choice of
  $\psi(p)$ is, for $0<p\le1$,
  \begin{equation}
    \psi(p)=\frac{1}{ l_\mathrm{max}} \;\ln\biggl(\frac{1}{p
      (d_\mathrm{max}-1)}\biggr).
  \end{equation}
  In particular, $g_{G,v}^{p}(n)$ decays exponentially if
  $p<1/(d_\mathrm{max}-1)$.
\end{lemma}

\begin{proof}
  A path with initial vertex $v$, which enters the complement of
  $B_n(v)$, contains at least $\widetilde{n}:= \lceil n/
  l_\mathrm{max}\rceil$ bonds, where $\lceil x\rceil$ is the smallest
  integer $\ge x$. With the notation in the proof of
  Lemma~\ref{combi}, the
  assertion now follows by noting that
  \begin{equation}
    g_{G,v}^{p}(n)\le \PP_G^p \bigl\{ \omega\in\Omega_{G}:
    W_{v}^{(\omega)}(\widetilde{n}) \ge 1 \bigr\} 
    \le \frac{d_{\mathrm{max}}}{d_\mathrm{max}-1}
    \;[p(d_{\mathrm{max}} -1)]^{\widetilde{n}}.  
  \end{equation}
\end{proof}

For a graph $G=(\cV,\cE)$ with uniformly discrete vertex set,
exponential decay of $g_{G,v}^{p}(n)$ implies that the mean cluster
size $\Chi_{G,v}(p) := \mathbb{E}_{G}^{p}\{|C_{v}|\}$ is finite, as
follows from the argument in \cite[p.~89]{Gri99}. In fact, this
argument yields an estimate which is uniform in $v\in\mathcal V$. We
state the result as

\begin{lemma}
  \label{finitechi}
  Let $G=(\cV,\cE)$ be a graph with uniformly discrete vertex set of
  radius $r>0$.  Assume that $G$ has maximal vertex degree $d_\mathrm{max}
  \in \mathbb{N} \setminus\{1\}$ and maximal edge length
  $l_\mathrm{max}<\infty$. Then, for every $p \in [0, \frac{1}{
    d_\mathrm{max}-1}[$, there exists a constant $\chi(p) \in
  ]0,\infty[$, which depends only on $r$, $\d_{\mathrm{max}}$ and
  $l_{\mathrm{max}}$ otherwise, such that the mean cluster size
  $\Chi_{G,v}(p) = \mathbb{E}_{G}^{p}\{|C_{v}|\}$ satisfies
  \begin{equation}
    \Chi_{G,v}(p)\le\chi(p)<\infty,
  \end{equation}
  uniformly in $v\in\cV$.\qed
\end{lemma}

Lemma~\ref{expsmallp} and Lemma~\ref{finitechi} 
can be used to prove exponential decay of the cluster size
distribution.

\begin{theorem}
  \label{theo:exp}
  Let $G=(\cV,\cE)$ be a graph with uniformly discrete vertex set of
  radius $r>0$. Assume that $G$ has maximal vertex degree
  $d_\mathrm{max} \in \mathbb{N} \setminus\{1\}$ and maximal edge
  length $l_\mathrm{max}<\infty$. Then, for every $p \in [0, \frac{1}{
    d_\mathrm{max}-1}[$, there exists a constant $\lambda(p) \in
  ]0,\infty[$, which depends only on $r$, $\d_{\mathrm{max}}$ and
  $l_{\mathrm{max}}$ otherwise, such that
  \begin{equation}
    {\mathbb P}_{G}^{p}\bigl(|C_v|\ge n \bigr)\le 2\, \e^{-n\lambda(p)}
  \end{equation}
  for all $n\in\mathbb N$, uniformly in $v\in\mathcal V$.
\end{theorem}

\begin{proof}[Sketch of the proof]
  Proceed along the lines of \cite[Thm.~6.75]{Gri99}. In the estimates,
  replace $\Chi_{G,v}(p)$ by its uniform bound $\Chi(p)$. The decay
  rate thus obtained is $\lambda(p) = [2\Chi(p)^{2}]^{-1}$.
\end{proof} 

The conclusion in the above theorem still holds if, instead of a
single graph $G$, we consider a set of graphs
$\mathcal G$ with the above properties and the associated dynamical
system $X_\mathcal{G}$. This setup is used in
Section~\ref{secLif}.

\begin{cor}
  \label{together}
  Let $\mathcal G$ be a set of graphs, whose vertex sets are uniformly
  discrete of radius $r>0$, which have maximum vertex degree
  $d_{\mathrm{max}}\in\mathbb{N}\setminus \{1\}$ and maximal edge
  length $l_\mathrm{max}<\infty$. Then, for every $p \in [0, \frac{1}{
    d_\mathrm{max}-1}[$ there exists a constant $\lambda(p) \in
  ]0,\infty[$ such that
  \begin{equation}
    {\mathbb P}_{G}^{p}\bigl(|C_v|\ge n \bigr)\le 2\, \e^{-n\lambda(p)}
  \end{equation}
  for all $n\in\mathbb{N}$,  uniformly in $G \in X_{\cG}$ and in
  $v\in\cV_{G}$.
  \qed
\end{cor}

\section*{Acknowledgements}

We are grateful to Michael Baake and Daniel Lenz for stimulating
discussions. We also thank Daniel Lenz and Ivan Veseli\'c for sending
us a version of their manuscript \cite{LeVe07} before making it
public. This work was supported by the German Research Council (DFG),
within the CRC 701.


\end{document}